\documentclass[11pt,leqno]{article}
\usepackage{a4}
\usepackage[T1]{fontenc}
\usepackage[english]{babel}
\usepackage{latexsym}
\usepackage{amsmath}
\usepackage{amssymb}
\usepackage{mathrsfs}
\usepackage{color}
\usepackage{cite}
\usepackage{titling}
\makeatletter
\@addtoreset{equation}{section} 
\makeatother
\definecolor{move}{rgb}{.3,.1,.8}
\newcommand{\ud}{\mathrm{d}}
\newcommand{\de}{\mathrm{D}}
\newcommand{\supp}{\mathrm{supp}}
\newcommand{\re}{\mathrm{Re}}

\newcommand{\e}{\mathrm{e}}
\newcommand{\id}{\mathrm{Id}}
\newcommand{\R}{\mathbb{R}}
\newcommand{\N}{\mathbb{N}}

\newcommand{\Ci}{\mathscr{C}^{\infty}}

\newenvironment{pr}{\vspace{5pt}\textbf{{\small Proof :}}\\}{\hspace{\stretch{1}}\rule{1ex}{1ex}\vspace{5pt}}
\newtheorem{thm}{Theorem}[section]
\newtheorem{pro}{Proposition}[section]

\newtheorem{lem}{Lemma}[section]

\usepackage{fancyhdr}
\fancypagestyle{plain}{\fancyhead{}}
\title{Energy decay estimates of elastic transmission wave/beam systems with a local Kelvin-Voigt damping}
\author{{FATHI HASSINE}\\\\  \textit{UR Analysis and Control of PDE 13ES64}\\ \textit{Department of Mathematics, Faculty of Sciences of Monastir}\\ \textit{University of Monastir, 5019 Monastir, Tunisia}\\ \textit{email:} \texttt{fathi.hassine@fsm.rnu.tn}}
\date{}
\begin{document}
\maketitle
\begin{center}
\abstract{
We consider a beam and a wave equations coupled on an elastic beam through transmission conditions. 
The damping which is locally distributed acts through one of the two equations only; its effect is transmitted to the other equation through the coupling. First we consider the case where the dissipation acts through the beam equation. Using a recent result of Borichev and Tomilov on polynomial decay characterization of bounded semigroups we provide a precise decay estimates showing that the energy of this coupled system decays polynomially as the time variable goes to infinity. Second, we discuss the case where the damping acts through the wave equation. Proceeding as in the first case, we prove that this system is also polynomially stable and we provide precise polynomial decay estimates for its energy. Finally, we show the lack of uniform exponential decay of solutions for both models.
}
\end{center}


\vspace{10pt}
\section{Introduction and motivation}
Consider a clamped elastic beam of length $L$. One segment of the beam is made of a viscoelastic material with Kelvin-Voigt constitutive relation. The longitudinal and transversal vibration of the beam can be described by the following equations
\begin{equation*}
\left\{\begin{array}{ll}
\rho\ddot{u}(x,t)-(pu'+d(x)\chi_{(\alpha,\beta)}\dot{u}')'(x,t)=0&\text{in }(0,L)\times(0,+\infty)
\\
u(0,t)=u(L,t)=0&\text{for }t\in(0,+\infty)
\\
u(x,0)=u_{1}(x),\;\dot{u}(x,0)=u_{2}(x)&\text{on }(0,L)
\end{array}\right.
\end{equation*}
and
\begin{equation*}
\left\{\begin{array}{ll}
\rho\ddot{w}(x,t)+(pw''+d(x)\chi_{(\alpha,\beta)}\dot{w}'')''(x,t)=0&\text{in }(0,L)\times(0,+\infty)
\\
w(0,t)=w'(0,t)=w(L,t)=w'(L,t)=0&\text{for }t\in(0,+\infty)
\\
w(x,0)=w_{0}(x),\;\dot{w}_{1}(x,0)=w_{1}(x)&\text{on }(0,L),
\end{array}\right.
\end{equation*}
where $u$ and $w$ represent respectively the longitudinal and transversal displacement of the beam in the interval $(0,L)$ (the prime denotes the space derivative and the dot denotes the time derivative). The coefficient functions $p,\,\rho$ and $d$ are strictly positive and in $L^{\infty}(0,L)$ and $0\leq\alpha<\beta\leq L$ with $\chi_{(\alpha,\beta)}$ being the characteristic function of the interval $(\alpha,\beta)$.

It is well known that the energy of the solutions of the system describing the longitudinal vibration of the beam is polynomially and not exponentially stable however the one of the transversal vibration of the beam is exponentially stable (see~\cite{LL} and~\cite{LR1}). One question of interest then is how the stability properties are affected if we couple the exponentially stable beam equations to the conservative wave equations and if we couple the polynomially stable wave equations to the conservative beam equations by transmission conditions. That is, we wonder how these properties are affected if we consider the two following systems
\begin{equation}\label{pwkv1}
\left\{\begin{array}{ll}
\ddot{u}_{1}(x,t)-(u_{1}'+D_{a}\dot{u}_{1}')'(x,t)=0&\text{in }(0,l)\times(0,+\infty)
\\
\ddot{u}_{2}(x,t)+u_{2}''''(x,t)=0&\text{in }(l,L)\times(0,+\infty)
\\
u_{1}(l,t)=u_{2}(l,t)&\text{for }t\in(0,+\infty)
\\
u_{2}'(l,t)=0&\text{for }t\in(0,+\infty)
\\
u_{1}'(l,t)+u_{2}'''(l,t)=0&\text{for }t\in(0,+\infty)
\\
u_{1}(0,t)=u_{2}(L,t)=u_{2}'(L,t)=0&\text{for }t\in(0,+\infty)
\\
u_{1}(x,0)=u_{1}^{0}(x),\;\dot{u}_{1}(x,0)=u_{1}^{1}(x)&\text{on }(0,l)
\\
u_{2}(x,0)=u_{2}^{0}(x),\;\dot{u}_{2}(x,0)=u_{2}^{1}(x)&\text{on }(l,L),
\end{array}\right.
\end{equation}
and
\begin{equation}\label{pwkv2}
\left\{\begin{array}{ll}
\ddot{w}_{1}(x,t)+(w_{1}''+D_{b}\dot{w}_{1}'')''(x,t)=0&\text{in }(0,l)\times(0,+\infty)
\\
\ddot{w}_{2}(x,t)-w_{2}''(x,t)=0&\text{in }(l,L)\times(0,+\infty)
\\
w_{1}(l,t)=w_{2}(l,t)&\text{for }t\in(0,+\infty)
\\
w_{1}'(l,t)=0&\text{for }t\in(0,+\infty)
\\
w_{1}'''(l,t)+w_{2}'(l,t)=0&\text{for }t\in(0,+\infty)
\\
w_{1}(0,t)=w_{1}'(0,t)=w_{2}(L,t)=0&\text{for }t\in(0,+\infty)
\\
w_{1}(x,0)=w_{1}^{0}(x),\;\dot{w}_{1}(x,0)=w_{1}^{1}(x)&\text{on }(0,l)
\\
w_{2}(x,0)=w_{2}^{0}(x),\;\dot{w}_{2}(x,0)=w_{2}^{1}(x)&\text{on }(l,L),
\end{array}\right.
\end{equation}
where $w_{1}$ and $u_{2}$ represent the transversal displacement of the beam and $w_{2}$ and $u_{1}$ represent the longitudinal one and $D_{a}=a(x)\chi_{(\alpha,\beta)}$ and $D_{b}=b(x)\chi_{(\alpha,\beta)}$ with $0<\alpha<\beta<l<L$ and $a(x),\,b(x)\geq c_{0}>0$ in $(\alpha,\beta)$.

The third, the fourth and the fifth equation of~\eqref{pwkv1} and~\eqref{pwkv2} are called the transmission conditions. The first one is known as the continuity transmission condition, the second is described by the fact that the slope of the beam is null (This can be hold for example by imaging the beam is along the $x$ axis and deflects in the $y$ direction. This can be shown for instance if a clamped end on a sliding bearing that slides in the $y$ direction or also on a clamp at the end of a speedometer cable\footnote{Speedometer cable is flexible in bending and stiff in torsion. Thus, when held off to the side it can allow the beam to bend without allowing it to rotate at the end.}) and the third, means that the two forces which are the shear force of the beam and the stress of the string are such that one cancels the other 
. Different transmission conditions have been treated in~\cite{AN} for the thin plate model, and for the longitudinal and transversal vibrations of the Euler-Bernoulli beam~\cite{Ha1}, more natural transmission conditions have been taken into account.

Now the questions of interest are, is the full above systems stable and, if so, at which rate? The energy of a solution of~\eqref{pwkv1} and~\eqref{pwkv2} at the time $t\geq 0$ are defined respectively by
\[
E_{1}(t)=\frac{1}{2}\left(\int_{0}^{l}\Big(|\dot{u}_{1}(x,t)|^{2}+|u_{1}'(x,t)|^{2}\Big)\,\ud x+\int_{l}^{L}\Big(|\dot{u}_{2}(x,t)|^{2}+|u_{2}''(x,t)|^{2}\Big)\,\ud x\right),
\]
and
\[
E_{2}(t)=\frac{1}{2}\left(\int_{0}^{l}\Big(|\dot{w}_{1}(x,t)|^{2}+|w_{1}''(x,t)|^{2}\Big)\,\ud x+\int_{l}^{L}\Big(|\dot{w}_{2}(x,t)|^{2}+|w_{2}'(x,t)|^{2}\Big)\,\ud x\right).
\]
By a formal calculation we can show that for all $t_{1},\,t_{2}>0$ we have
\[
E_{1}(t_{2})-E_{1}(t_{1})=-\int_{t_{1}}^{t_{2}}\!\!\!\int_{\alpha}^{\beta}D_{a}|\dot{u}_{1}'(x,t)|^{2}\,\ud x\,\ud t,
\]
and
\[
E_{2}(t_{2})-E_{2}(t_{1})=-\int_{t_{1}}^{t_{2}}\!\!\!\int_{\alpha}^{\beta}D_{b}|\dot{w}_{1}''(x,t)|^{2}\,\ud x\,\ud t.
\]
This mean that the energy is decreasing over the time for both systems.

This last years the study of the stabilization problem for coupled systems has attracted a lot of attention e.g. \cite{Al}, \cite{ARSV}, \cite{Ava}, \cite{AvLa}, \cite{BEPS}, \cite{BRA}, \cite{AN}, \cite{AV}, \cite{Ha1}, \cite{Ha3}, \cite{Ha4}, \cite{D}, \cite{I}, \cite{LZ}, \cite{RZZ}, \cite{T}, \cite{ZZ1}, \cite{ZZ2} and \cite{ZZ3}. The systems discussed in those paper involve thermoelastic systems, fluid-structure interaction systems, and coupled wave-wave, plate-plate, or plate-wave equations. The techniques developed for a such coupled systems are very diverse. We can cite the approach based on the use of a Riesz basis, the frequency method based on Carleman estimates or on the multiplier method, the observability inequality, the approach based on spectral analysis\dots

What makes those kind of problems interesting is that the damping acts through only one equation. In addition, in our cases we are dealing with a locally distributed damping. This leads to technical difficulties when one tries to estimate the energy of the undamped equation. Our main purpose in this work is to develop a device that will help us to estimate the decay rate of the energy. Due to the locally distributed and unbounded nature of the damping, we use a frequency domain method and combine a contradiction argument with the multiplier technique to carry out a special analysis for the resolvent. Especially we will show that the energy of solutions of~\eqref{pwkv1} and~\eqref{pwkv2} is polynomially stable and due to the presence of the wave equation it is not exponentially stable in both cases.

This paper is organized as follows. In section~\ref{pwkv122} we gives the main results. In section~\ref{pwkv94} we  discuss the case where the damping acts through the plate equation. In section~\ref{pwkv95} is devoted for the case where the damping acts through the wave equation. In section~\ref{pwkv119} we show that the solutions of systems~\eqref{pwkv1} and~\eqref{pwkv2} are not exponentially stable.
\section{Preliminary and main results}\label{pwkv122}
Let $H=L^{2}(0,l)\times L^{2}(l,L)$ with the norm
$$\|v\|_{H}=\|(v_{1},v_{2})\|_{H}=\left(\int_{0}^{l}|v_{1}|^{2}(x)\,\ud x+\int_{l}^{L}|v_{2}|^{2}(x)\,\ud x\right)^{\frac{1}{2}}$$
and 
\begin{equation*}
\begin{split}
V_{1}=\{w=(w_{1},w_{2})\in H^{2}(0,l)\times H^{1}(l,L): w_{1}(l)=w_{2}(l),\,w_{1}(0)=w_{2}(L)=w_{1}'(0)=0\}
\end{split}
\end{equation*}
with the norm
$$\|w\|_{V_{1}}=\|(w_{1},w_{2})\|_{V_{1}}=\left(\int_{0}^{l}|w_{1}''|^{2}(x)\,\ud x+\int_{l}^{L}|w_{2}'|^{2}(x)\,\ud x\right)^{\frac{1}{2}}.$$
Define $\mathcal{H}_{1}=V_{1}\times H$ with the norm $\|(w,v)\|_{\mathcal{H}_{1}}^{2}=\|w\|_{V_{1}}^{2}+\|v\|_{H}^{2}$. Then $\mathcal{H}_{1}$ is a Hilbert space and we define
\begin{equation*}
\begin{split}
\mathcal{D}(\mathcal{A}_{1})=\{(w,v)=((w_{1},w_{2}),(v_{1},v_{2}))\in\mathcal{H}_{1}:v\in V_{1},\,w_{1}''+D_{b}.v_{1}''\in H^{2}(0,l),
\\
w_{2}'\in H^{1}(l,L),\,w_{1}'''(l)+w_{2}'(l)=0\}
\end{split}
\end{equation*}
and
$$\mathcal{A}_{1}(w,v)=\mathcal{A}_{1}((w_{1},w_{2}),(v_{1},v_{2}))=(v_{1},v_{2},-(w_{1}''+D_{a}.v_{1}'')'',w_{2}'').$$
Thus,~\eqref{pwkv2} can be rewritten as an abstract evolution equation on $\mathcal{H}_{1}$,
$$(\dot{w}(t),\dot{v}(t))=\mathcal{A}_{1}(w(t),v(t)),\quad (w_{1}(0),w_{2}(0),v_{1}(0),v_{2}(0))=(w_{1}^{0},w_{2}^{0},w_{1}^{1},w_{2}^{1}).$$
\begin{pro}\label{pwkv96}
The linear operator $\mathcal{A}_{1}$ generates a $C_{0}$-semigroup of contractions $e^{t\mathcal{A}_{1}}$ on $\mathcal{H}_{1}$, in particular there exists a unique solution of~\eqref{pwkv2} which can be expressed by means of a semigroup on $\mathcal{H}_{1}$ having the following regularity of the solution
$$
\left(\begin{array}{l}
w_{1}
\\
w_{2}
\\
\dot{w}_{1}
\\
\dot{w}_{2}
\end{array}\right)
\in C([0,+\infty[,\mathcal{D}(\mathcal{A}_{1}))\cap C^{1}([0,+\infty[,\mathcal{H}_{1}),$$
if $(w_{1}^{0},w_{2}^{0},w_{1}^{1},w_{2}^{1})\in \mathcal{D}(\mathcal{A}_{1})$ and a mild solution
$$
\left(\begin{array}{l}
w_{1}
\\
w_{2}
\\
\dot{w}_{1}
\\
\dot{w}_{2}
\end{array}\right)
\in C([0,+\infty[,\mathcal{H}_{1}),$$
if $(w_{1}^{0},w_{2}^{0},w_{1}^{1},w_{2}^{1})\in\mathcal{H}_{1}$. Further the semigroup $e^{\mathcal{A}_{1}t}$ is strongly stable i.e
$$
\lim_{t\longrightarrow+\infty}\left\|e^{\mathcal{A}_{1}t}\left(\begin{array}{l}
w_{1}^{0}
\\
w_{2}^{0}
\\
w_{1}^{1}
\\
w_{2}^{1}
\end{array}\right)\right\|_{\mathcal{H}_{1}}=0\qquad\forall\,\left(\begin{array}{l}
w_{1}^{0}
\\
w_{2}^{0}
\\
w_{1}^{1}
\\
w_{2}^{1}
\end{array}\right)\in\mathcal{H}_{1}.
$$
\end{pro}
\begin{pr}
The operator $\mathcal{A}_{1}$ is dissipative by the fact,
$$
\re\left\langle\mathcal{A}_{1}\left(\begin{array}{l}
w_{1}
\\
w_{2}
\\
v_{1}
\\
v_{2}
\end{array}\right),\left(\begin{array}{l}
w_{1}
\\
w_{2}
\\
v_{1}
\\
v_{2}
\end{array}\right)\right\rangle_{\mathcal{H}_{1}}=-\int_{0}^{l}b|v_{1}''|\,\ud x\leq 0.
$$
Moreover, by Lax-Miligram Theorem $(\id-\mathcal{A}_{1})$ is onto. Consequently, $\mathcal{A}_{1}$ generates a $C_{0}$-semigroup of contractions $e^{t\mathcal{A}_{1}}$ on $\mathcal{H}_{1}$ (see~\cite{P}).

For the the strong stability, it is easy to show that there is no point spectrum of $\mathcal{A}_{1}$ on the imaginary axis, i.e. $i\R\cap\sigma_{p}(\mathcal{A}_{1})\neq\emptyset$. Then it follow from~\cite[Lemma 4.1]{CLL} that the resolvent set of $\mathcal{A}_{1}$ contains the imaginary axis. Thus the result follow easily from~\cite{Ben}.
\end{pr}

Our first main result is now given by the following
\begin{thm}\label{pwkv120}
The semigroup $e^{\mathcal{A}_{1}t}$ is polynomially stable and in particular, there exist $M>0$ such that for all $t\geq 0$
\[
\|e^{\mathcal{A}_{1}t}(w_{1}^{0},w_{2}^{0},w_{1}^{1},w_{2}^{1})\|_{\mathcal{H}_{1}}\leq\frac{M}{1+t}\|(w_{1}^{0},w_{2}^{0},w_{1}^{1},w_{2}^{1})\|_{\mathcal{D}(\mathcal{A}_{1})},\,\forall\,(w_{1}^{0},w_{2}^{0},w_{1}^{1},w_{2}^{1})\in\mathcal{D}(\mathcal{A}_{1}).
\]
Besides the semigroup $e^{t\mathcal{A}_{1}}$ is not exponentially stable.
\end{thm}

We focus now to the system~\eqref{pwkv1} and we define $V_{2}$ by
\begin{equation*}
\begin{split}
V_{2}=\{u=(u_{1},u_{2})\in H^{1}(0,l)\times H^{2}(l,L):u_{1}(l)=u_{2}(l),\,u_{2}'(l)=0,
\\
u_{1}(0)=0,\,u_{2}(L)=0,\,u_{2}'(L)=0\}
\end{split}
\end{equation*}
with the norm
$$\|u\|_{V_{2}}=\|(u_{1},u_{2})\|_{V_{2}}=\left(\int_{0}^{l}|u_{1}'|^{2}(x)\,\ud x+\int_{l}^{L}|u_{2}''|^{2}(x)\,\ud x\right)^{\frac{1}{2}}.$$
We define $\mathcal{H}_{2}=V_{2}\times H$ with norm $\|(u,v)\|_{\mathcal{H}_{2}}^{2}=\|u\|_{V_{2}}^{2}+\|v\|_{H}^{2}$. Then $\mathcal{H}_{2}$ is a Hilbert space in which we define the operator $\mathcal{A}_{2}$ by
$$\mathcal{A}_{2}(u,v)=\mathcal{A}_{2}((u_{1},u_{2}),(v_{1},v_{2}))=(v_{1},v_{2},(u'_{1}+D_{a}v_{1}')',-u_{2}'''')$$
with domain
\begin{equation*}
\begin{split}
\mathcal{D}(\mathcal{A}_{2})=\{(u,v)=((u_{1},u_{2}),(v_{1},v_{2}))\in\mathcal{H}_{2}: v\in V_{2},\,u_{1}'+D_{a}v_{1}\in H^{1}(0,l),\,u_{2}''\in H^{2}(l,L),
\\
u_{1}'(l)+u_{2}'''(l)=0\}.
\end{split}
\end{equation*}
Then~\eqref{pwkv1} can be rewritten as an abstract evolution equation on $\mathcal{H}_{2}$,
$$(\dot{u}(t),\dot{v}(t))=\mathcal{A}_{2}(u(t),v(t)),\quad (u_{1}(0),u_{2}(0),v_{1}(0),v_{2}(0))=(u_{1}^{0},u_{2}^{0},u_{1}^{1},u_{2}^{1}).$$
Similar to Proposition~\ref{pwkv96} we can prove here also the following
\begin{pro}
The linear operator $\mathcal{A}_{2}$ generates a $C_{0}$-semigroup of contractions $e^{t\mathcal{A}_{2}}$ on $\mathcal{H}_{2}$, in particular there exists a unique solution of~\eqref{pwkv1} which can be expressed by means of a semigroup on $\mathcal{H}_{2}$ having the following regularity of the solution
$$
\left(\begin{array}{l}
u_{1}
\\
u_{2}
\\
\dot{u}_{1}
\\
\dot{u}_{2}
\end{array}\right)
\in C([0,+\infty[,\mathcal{D}(\mathcal{A}_{2}))\cap C^{1}([0,+\infty[,\mathcal{H}_{2}),$$
if $(u_{1}^{0},u_{2}^{0},u_{1}^{1},u_{2}^{1})\in \mathcal{D}(\mathcal{A}_{2})$ and a mild solution
$$
\left(\begin{array}{l}
u_{1}
\\
u_{2}
\\
\dot{u}_{1}
\\
\dot{u}_{2}
\end{array}\right)
\in C([0,+\infty[,\mathcal{H}_{2}),$$
if $(u_{1}^{0},u_{2}^{0},u_{1}^{1},u_{2}^{1})\in\mathcal{H}_{2}$. Further, the semigroup $e^{\mathcal{A}_{2}t}$ is strongly stable i.e
$$
\lim_{t\longrightarrow+\infty}\left\|e^{\mathcal{A}_{2}t}\left(\begin{array}{l}
u_{1}^{0}
\\
u_{2}^{0}
\\
u_{1}^{1}
\\
u_{2}^{1}
\end{array}\right)\right\|_{\mathcal{H}_{2}}=0\qquad\forall\,\left(\begin{array}{l}
u_{1}^{0}
\\
u_{2}^{0}
\\
u_{1}^{1}
\\
u_{2}^{1}
\end{array}\right)\in\mathcal{H}_{2}.
$$
\end{pro}

Assuming now that the function $a\in\mathrm{C}(\alpha,\beta)$ and we give now the second main result 
\begin{thm}\label{pwkv121}
Under the above assumptions on the coefficients of~\eqref{pwkv1} the semigroup $e^{\mathcal{A}_{2}t}$ is polynomially stable and in particular, there exist $M>0$ such that
\[
\|e^{\mathcal{A}_{2}t}(u_{1}^{0},u_{2}^{0},u_{1}^{1},u_{2}^{1})\|_{\mathcal{H}_{2}}\leq\frac{M}{1+t}\|(u_{1}^{0},u_{2}^{0},u_{1}^{1},u_{2}^{1})\|_{\mathcal{D}(\mathcal{A}_{2})},\quad\forall\,(u_{1}^{0},u_{2}^{0},u_{1}^{1},u_{2}^{1})\in\mathcal{D}(\mathcal{A}_{2}),\,\forall\,t\geq 0.
\]
Besides the semigroup $e^{t\mathcal{A}_{2}}$ is not exponentially stable.
\end{thm}
\section{Damping arising from the transversal motion}\label{pwkv94}
The purpose of this section is to prove the fist part of Theorem~\ref{pwkv120}. We need only to verify the condition for a semigroup of contractions on a Hilbert space being polynomially stable (see~\cite{BT}), i.e.,
\begin{equation}\label{pwkv4}
\sup\{|\lambda|^{-1}.\|(i\lambda-\mathcal{A}_{1})^{-1}\|;\quad\lambda\in\R\}<+\infty.
\end{equation}
Suppose that~\eqref{pwkv4} is not true. By the continuity of the resolvent and the resonance theorem, there exist $\lambda_{n}\in\R$, $(w_{n},v_{n})=((w_{1,n},w_{2,n}),(v_{1,n},v_{2,n}))\in\mathcal{D}(\mathcal{A}_{1})$, for all $n\in\N$ such that
\begin{equation}\label{pwkv3}
\|(w_{n},v_{n})\|_{\mathcal{H}_{1}}=1,\;\lim_{n\rightarrow+\infty}|\lambda_{n}|=+\infty
\end{equation}
and
\begin{equation}\label{pwkv15}
|\lambda_{n}|.(i\lambda_{n}-\mathcal{A}_{1})(w_{n},v_{n})=(f_{n},g_{n})=((f_{1,n},f_{2,n}),(g_{1,n},g_{2,n}))\,\longrightarrow\, 0\text{ in }\mathcal{H}_{1}.
\end{equation}
This implies
\begin{eqnarray}
|\lambda_{n}|.(i\lambda_{n}w_{1,n}-v_{1,n})=f_{1,n}\,\longrightarrow\,0\text{ in }H^{2}(0,l),\label{pwkv5}
\\
|\lambda_{n}|.(i\lambda_{n}w_{2,n}-v_{2,n})=f_{2,n}\,\longrightarrow\,0\text{ in }H^{1}(l,L),\label{pwkv6}
\\
|\lambda_{n}|.(i\lambda_{n}.v_{1,n}-M_{n}'')=g_{1,n}\,\longrightarrow\,0\text{ in }L^{2}(0,l),\label{pwkv7}
\\
|\lambda_{n}|.(i\lambda_{n}.v_{2,n}-w_{2,n}'')=g_{2,n}\,\longrightarrow\,0\text{ in }L^{2}(l,L),\label{pwkv8}
\end{eqnarray}
where 
$$
M_{n}=-(w_{1,n}''+D_{b}v_{1,n}'').
$$
For $0\leq x\leq l$ define
\begin{equation}\label{pwkv9}
J(\psi)(x)=\int_{x}^{l}\!\!\int_{s}^{l}\psi(\tau)\,\ud\tau\,\ud s
\end{equation}
and
\begin{equation}\label{pwkv10}
y_{n}=\frac{1}{i\lambda_{n}}[|\lambda_{n}|^{\frac{1}{2}}.M_{n}+|\lambda_{n}|^{-\frac{1}{2}}.J(g_{1,n})].
\end{equation}
Comparing~\eqref{pwkv9} and~\eqref{pwkv10} we have
\begin{equation}\label{pwkv11}
y_{n}''=|\lambda_{n}|^{\frac{1}{2}}.v_{1,n}.
\end{equation}
The rest of the proof depends on the following two lemmas. Let $\omega_{n}=\sqrt{|\lambda_{n}|}$.
\begin{lem}
The function $y_{n}$ defined above has the following properties:
\begin{align}
y_{n}\,\longrightarrow\,0\quad\text{ in }\; H^{4}(\alpha,\beta),\label{pwkv12}
\\
\lambda_{n}.y_{n}\,\longrightarrow\,0\quad\text{ in }\; L^{2}(\alpha,\beta),\label{pwkv13}
\\
\omega_{n}.y_{n}\,\longrightarrow\,0\quad\text{ in }\; H^{2}(\alpha,\beta).\label{pwkv14}
\end{align}
\end{lem}
\begin{pr}
From~\eqref{pwkv15}, we have
\begin{equation}\label{pwkv16}
\re\langle|\lambda_{n}|.(i\lambda_{n}-\mathcal{A}_{1})(w_{n},v_{n}),(w_{n},v_{n})\rangle_{\mathcal{H}}=|\lambda_{n}|.\int_{\alpha}^{\beta}b(x)|v_{1,n}''|^{2}\,\ud x\,\longrightarrow\,0.
\end{equation}
Therefore, from~\eqref{pwkv5} we have
\begin{equation}\label{pwkv17}
\omega_{n}.M_{n}\,\longrightarrow\,0\quad\text{ in }\;L^{2}(\alpha,\beta)
\end{equation}
and
\begin{equation}\label{pwkv18}
\frac{1}{|\lambda_{n}|}\|\psi.v_{n}\|_{V_{1}}=O(1)
\end{equation}
for every $\psi\in\Ci([0,l])$, such that $\supp(\psi)\subset[0,l]$.

Equations~\eqref{pwkv7},~\eqref{pwkv17} and~\eqref{pwkv18} imply that
\begin{equation}\label{pwkv19}
|\lambda_{n}|.\int_{\alpha}^{\beta}\psi.|v_{1,n}|^{2}\,\ud x\,\longrightarrow\,0,\quad\forall\,\psi\in\Ci([0,l]),\;\supp(\psi)\subset(\alpha,\beta).
\end{equation}
Applying the interpolation theorem involving compact subdomain~\cite[Theorem 4.23]{A} we find that~\eqref{pwkv16} and~\eqref{pwkv19} imply
\begin{equation}\label{pwkv20}
\omega_{n}.v_{1,n}\,\longrightarrow\,0\quad\text{ in }\;H^{2}(\alpha,\beta).
\end{equation}
Thus,~\eqref{pwkv11} yields
\begin{equation}\label{pwkv21}
\int_{\alpha}^{\beta}|y_{n}''''|^{2}\,\ud x\,\longrightarrow\,0.
\end{equation}
On the other hand,~\eqref{pwkv13} follow from
\begin{equation}\label{pwkv22}
i\lambda_{n}y_{n}=\omega_{n}^{-1}.J(g_{1,n})-\frac{1+i\lambda_{n}D_{b}}{i\lambda_{n}}\omega_{n}.v_{n}''-\frac{1}{i|\lambda_{n}|.\lambda_{n}}f_{1,n}''\,\longrightarrow\,0\quad\text{ in }\;L^{2}(\alpha,\beta).
\end{equation}
Since $|\lambda_{n}|\,\longrightarrow\,+\infty$, we obtain that $y_{n}\,\longrightarrow\,0$ in $L^{2}(\alpha,\beta)$. This combined with~\eqref{pwkv21} yields~\eqref{pwkv12}. From the interpolation inequality~\cite[Theorem 4.17]{A}, we also have~\eqref{pwkv14}.
\end{pr}
\begin{lem}\label{pwkv111}
The functions $w_{1,n}\in H^{4}(0,\alpha)\cap H^{4}(\beta,l)$, for all $n\in\N$ have the following properties:
\begin{align}
\omega_{n}^{3}(|w_{1,n}(\alpha)|+|w_{1,n}'(\alpha)|+|w_{1,n}(\beta)|+|w_{1,n}'(\beta)|)\,\longrightarrow\,0,\label{pwkv23}
\\
|w_{1,n}'''(\alpha^{-})|\,\longrightarrow\,0\quad\text{ and }\quad|w_{1,n}'''(\beta^{+})|\,\longrightarrow\,0,\label{pwkv25}
\\
\alpha|w_{1,n}''(\alpha^{-})|^{2}+(l-\beta)|w_{1,n}''(\beta^{+})|^{2}+\int_{l}^{L}|\lambda_{n}.w_{2,n}|^{2}+3\int_{l}^{L}|w_{2,n}'|^{2}\,\ud x\longrightarrow\,2,\label{pwkv103}
\\
\alpha|w_{1,n}''(\alpha^{-})|^{2}+(l-\beta)|w_{1,n}''(\beta^{+})|^{2}+(L-l)(|w_{2,n}'(l)|^{2}+|\lambda_{n}.w_{2,n}(l)|^{2})+2\int_{l}^{L}\!\!\!|w_{2,n}'|^{2}\,\ud x\longrightarrow\,2,\label{pwkv42}
\\
\alpha|w_{1,n}''(\alpha^{-})|^{2}+(L-\beta)|w_{1,n}''(\beta^{+})|^{2}+(L-l)(|w_{2,n}'(l)|^{2}-|w_{1,n}''(l)|^{2})+2\int_{l}^{L}\!\!\!|w_{2,n}'|^{2}\,\ud x\longrightarrow\,2,\label{pwkv24}
\end{align}
and in particular we have
\begin{equation}\label{pwkv102}
|w_{1,n}''(\beta^{+})|^{2}-|w_{1,n}''(l)|^{2}-|\lambda_{n}.w_{1,n}(l)|^{2}\,\longrightarrow\,0.
\end{equation}
\end{lem}
\begin{pr}
Since $w_{1,n},\,v_{1,n}\in H^{2}(0,l)$, Sobolev's embedding theorem implies that they are also in $C^{1}(0,l)$. By~\eqref{pwkv5} and~\eqref{pwkv20} we have
\begin{equation}\label{pwkv26}
\omega_{n}^{3}.w_{1,n}\,\longrightarrow\,0\quad\text{ in }\;H^{2}(\alpha,\beta).
\end{equation}
Thus, $\omega_{n}^{3}.w_{1,n}$ converges to zero in $C^{1}([\alpha,\beta])$, which immediately leads to~\eqref{pwkv23}.

Note that $M_{n}=-w_{1,n}''$ on $(0,\alpha)\cup(\beta,l)$. From the definition of the domain of $\mathcal{A}_{1}$, we know that $w_{1,n}\in H^{4}(0,\alpha)$ and $w_{1,n}\in H^{4}(\beta,l)$. It follows from~\eqref{pwkv10} that
\begin{align}
\omega_{n}.w_{1,n}''(\alpha^{-})=(\omega_{n}^{-1}.J(g_{1,n})-i\lambda_{n}y_{n})(\alpha),\quad \omega_{n}.w_{1,n}''(\beta^{+})=(\omega_{n}^{-1}.J(g_{1,n})-i\lambda_{n}y_{n})(\beta),\label{pwkv27}
\\
\omega_{n}.w_{1,n}'''(\alpha^{-})=(\omega_{n}^{-1}.J(g_{1,n})-i\lambda_{n}y_{n})'(\alpha),\quad \omega_{n}.w_{1,n}'''(\beta^{+})=(\omega_{n}^{-1}.J(g_{1,n})-i\lambda_{n}y_{n})'(\beta).\label{pwkv28}
\end{align}
Dividing~\eqref{pwkv28} by $\omega_{n}$ we obtain~\eqref{pwkv25} by using~\eqref{pwkv14} in the previous lemma.

In ordre to prove~\eqref{pwkv103}-\eqref{pwkv24}, we substitute~\eqref{pwkv5} into~\eqref{pwkv7} and~\eqref{pwkv6} into~\eqref{pwkv8} to get
\begin{equation}\label{pwkv29}
\left\{\begin{array}{ll}
-\lambda_{n}^{2}.w_{1,n}-M_{n}''=|\lambda_{n}|^{-1}.(g_{1,n}+i\lambda_{n}.f_{1,n})&\text{in }(0,l)
\\
-\lambda_{n}^{2}.w_{2,n}-w_{2,n}''=|\lambda_{n}|^{-1}.(g_{2,n}+i\lambda_{n}.f_{2,n})&\text{in }(l,L).
\end{array}\right.
\end{equation}
We multiply the above equations by $\overline{w}_{1,n}$ and $\overline{w}_{2,n}$ respectively, then integrate by parts on $(0,L)$. This leads to
\begin{equation}\label{pwkv30}
\|\lambda_{n}.w_{n}\|_{H}^{2}-\|w_{n}\|_{V_{1}}^{2}\,\longrightarrow\,0.
\end{equation}
Here we have used~\eqref{pwkv3},~\eqref{pwkv5},~\eqref{pwkv6},~\eqref{pwkv7},~\eqref{pwkv8} and~\eqref{pwkv17}. Since $\|w_{n}\|_{V_{1}}^{2}+\|v_{n}\|_{H}^{2}=1$ and $i\lambda_{n}.w_{n}-v_{n}$ also converges to zero in $L^{2}(0,L)$. Hence~\eqref{pwkv30} implies that both $\|\lambda_{n}.w_{n}\|_{H}^{2}$ and $\|w_{n}\|_{V_{1}}^{2}$ must converge to $\displaystyle\frac{1}{2}$ as $n\,\longrightarrow\,+\infty$. This further leads to
\begin{equation}\label{pwkv31}
\begin{split}
\lim_{n\rightarrow+\infty}\int_{0}^{\alpha}|\lambda_{n}.w_{1,n}|^{2}\,\ud x+\int_{\beta}^{l}|\lambda_{n}.w_{1,n}|^{2}\,\ud x+\int_{l}^{L}|\lambda_{n}.w_{2,n}|^{2}\,\ud x=
\\
\lim_{n\rightarrow+\infty}\int_{0}^{\alpha}|w_{1,n}''|^{2}\,\ud x+\int_{\beta}^{l}|w_{1,n}''|^{2}\,\ud x+\int_{l}^{L}|w_{2,n}'|^{2}\,\ud x=\frac{1}{2}
\end{split}
\end{equation}
when~\eqref{pwkv26} is taken into account.

On the intervals $(0,\alpha)$, $(\beta,l)$ and $(l,L)$,~\eqref{pwkv29} becomes
\begin{equation}\label{pwkv32}
\left\{\begin{array}{ll}
-\lambda_{n}^{2}.w_{1,n}+w_{1,n}''''=|\lambda_{n}|^{-1}.(g_{1,n}+i\lambda_{n}.f_{1,n})&\text{in }(0,\alpha)\cup(\beta,l)
\\
-\lambda_{n}^{2}.w_{2,n}-w_{2,n}''=|\lambda_{n}|^{-1}.(g_{2,n}+i\lambda_{n}.f_{2,n})&\text{in }(l,L).
\end{array}\right.
\end{equation}
We multiply the first equation of~\eqref{pwkv32} respectively by $x\overline{w'}_{1,n}$ and $(L-x)\overline{w'}_{1,n}$ and integrate on $(0,\alpha)$ and $(\beta,l)$ respectively. Hence,
\begin{equation}\label{pwkv33}
-\int_{0}^{\alpha}\lambda_{n}^{2}w_{1,n}x\overline{w'}_{1,n}\,\ud x+\int_{0}^{\alpha}w_{1,n}''''x\overline{w'}_{1,n}\,\ud x=|\lambda_{n}|^{-1}.\int_{0}^{\alpha}(g_{1,n}+i\lambda_{n}f_{1,n})x\overline{w'}_{1,n}\,\ud x
\end{equation}
and
\begin{equation}\label{pwkv34}
\begin{split}
-\int_{\beta}^{l}\lambda_{n}^{2}w_{1,n}(L-x)\overline{w'}_{1,n}\,\ud x+\int_{\beta}^{l}w_{1,n}''''(L-x)\overline{w'}_{1,n}\,\ud x
\\
=|\lambda_{n}|^{-1}\int_{\beta}^{l}(g_{1,n}+i\lambda_{n}f_{1,n})(L-x)\overline{w'}_{1,n}\,\ud x.
\end{split}
\end{equation}
It is easy to see that the terms on the right hand side of~\eqref{pwkv33} and~\eqref{pwkv34} converges to zero. After a straightforward calculation (integration by parts), the two terms on the left hand side of~\eqref{pwkv33} and~\eqref{pwkv34} become
\begin{equation}\label{pwkv35}
-\re\left(\int_{0}^{\alpha}\lambda_{n}^{2}w_{1,n}x\overline{w'}_{1,n}\,\ud x\right)=-\omega_{n}^{4}\frac{\alpha}{2}|w_{1,n}(\alpha)|^{2}+\frac{1}{2}\int_{0}^{\alpha}|\lambda_{n}w_{1,n}|^{2}\,\ud x,
\end{equation}
\begin{equation}\label{pwkv36}
\begin{split}
\re\left(\int_{0}^{\alpha}w_{1,n}''''x\overline{w'}_{1,n}\,\ud x\right)=\re[(\alpha w_{1,n}'''(\alpha^{-})-w_{1,n}''(\alpha^{-}))\overline{w'}_{1,n}(\alpha)]
\\
+\frac{3}{2}\int_{0}^{\alpha}|w_{1,n}''|^{2}\,\ud x-\frac{\alpha}{2}|w_{1,n}''(\alpha^{-})|^{2},
\end{split}
\end{equation}
and
\begin{equation}\label{pwkv37}
\begin{split}
-\re\left(\int_{\beta}^{l}\lambda_{n}^{2}w_{1,n}(L-x)\overline{w'}_{1,n}\,\ud x\right)=-\frac{1}{2}\int_{\beta}^{l}|\lambda_{n}w_{1,n}|^{2}\,\ud x
\\
+\frac{\omega_{n}^{4}}{2}[(L-\beta)|w_{1,n}(\beta)|^{2}-(L-l)|w_{1,n}(l)|^{2}],
\end{split}
\end{equation}
\begin{equation}\label{pwkv38}
\begin{split}
\re\left(\int_{\beta}^{l}w_{1,n}''''(L-x)\overline{w'}_{1,n}\,\ud x\right)=\re[((L-l) w_{1,n}'''(l)+w_{1,n}''(l))\overline{w'}_{1,n}(l)]
\\
-\re[((L-\beta) w_{1,n}'''(\beta^{+})+w_{1,n}''(\beta^{+}))\overline{w'}_{1,n}(\beta)]-\frac{3}{2}\int_{\beta}^{l}|w_{1,n}''|^{2}\,\ud x
\\
+\frac{(L-\beta)}{2}|w_{1,n}''(\beta^{+})|^{2}-\frac{(L-l)}{2}|w_{1,n}''(l)|^{2}.
\end{split}
\end{equation}
After substituting these terms into the real part of~\eqref{pwkv33} and~\eqref{pwkv34} and applying~\eqref{pwkv23},~\eqref{pwkv25} and~\eqref{pwkv26}, we obtain
\begin{equation}\label{pwkv39}
\frac{1}{2}\int_{0}^{\alpha}|\lambda_{n}w_{1,n}|^{2}\,\ud x+\frac{3}{2}\int_{0}^{\alpha}|w_{1,n}''|^{2}\,\ud x-\frac{\alpha}{2}|w_{1,n}''(\alpha^{-})|^{2}\,\longrightarrow\,0,
\end{equation}
and
\begin{equation}\label{pwkv40}
\begin{split}
\frac{1}{2}\int_{\beta}^{l}|\lambda_{n}w_{1,n}|^{2}\,\ud x+\frac{3}{2}\int_{\beta}^{l}|w_{1,n}''|^{2}\,\ud x-\re[((L-l)w_{1,n}'''(l)+w_{1,n}''(l))\overline{w'}_{1,n}(l)]
\\
-\frac{(L-\beta)}{2}|w_{1,n}''(\beta^{+})|^{2}+\frac{(L-l)}{2}|w_{1,n}''(l)|^{2}+\frac{(L-l)}{2}\omega_{n}^{4}|w_{1,n}(l)|^{2}\,\longrightarrow\,0.
\end{split}
\end{equation}
Similarly, we can multiply the second equation of~\eqref{pwkv32} by $(L-x)\overline{w'}_{2,n}$ and integrate on $(l,L)$ to get
\begin{equation}\label{pwkv41}
\frac{1}{2}\int_{l}^{L}|\lambda_{n}w_{2,n}|^{2}\,\ud x+\frac{1}{2}\int_{l}^{L}|w_{2,n}'|^{2}\,\ud x-\frac{(L-l)}{2}\left(\omega_{n}^{4}|w_{2,n}(l)|^{2}+|w_{2,n}'(l)|^{2}\right)\,\longrightarrow\,0.
\end{equation}
Then~\eqref{pwkv24} follow by summing~\eqref{pwkv39},~\eqref{pwkv40} and~\eqref{pwkv41} and using~\eqref{pwkv31} and the transmission condition $w_{1,n}(l)=w_{2,n}(l)$.

Similarly in ordre to prove~\eqref{pwkv103} and~\eqref{pwkv42} we multiply the first equation of~\eqref{pwkv32} by $(l-x)\overline{w'}_{1,n}$ and integrating over $(\beta,l)$ then we obtain
\begin{equation}\label{pwkv99}
\frac{1}{2}\int_{\beta}^{l}|\lambda_{n}w_{1,n}|^{2}\ud x+\frac{3}{2}\int_{\beta}^{l}|w_{1,n}''|^{2}\ud x-\frac{(l-\beta)}{2}|w_{1,n}''(\beta^{+})|\,\longrightarrow\,0,
\end{equation}
hence~\eqref{pwkv103} follow by summing~\eqref{pwkv39} and~\eqref{pwkv99} and using~\eqref{pwkv31} when~\eqref{pwkv42} follow by summing~\eqref{pwkv39},~\eqref{pwkv41} and~\eqref{pwkv99} and using again~\eqref{pwkv31}.
Finally,~\eqref{pwkv102} follow easily by taking the difference of~\eqref{pwkv24} and~\eqref{pwkv42}.
\end{pr}

We will show now that
\begin{equation}\label{pwkv43}
|\omega_{n}w_{1,n}''(\alpha^{-})|^{2}\,\longrightarrow\,0\quad\text{ and }\quad|\omega_{n}w_{1,n}''(\beta^{+})|^{2}\,\longrightarrow\,0.
\end{equation}
Denote by
\begin{equation*}
\begin{array}{lll}
\displaystyle \widetilde{w}_{1,n}=\omega_{n}.w_{1,n}&\text{and}&\displaystyle F_{1,n}=\omega_{n}^{-1}(g_{1,n}+i\lambda_{n} f_{1,n}).
\end{array}
\end{equation*}
Then the first equation of~\eqref{pwkv32} can be written as
\begin{equation}\label{pwkv44}
\begin{array}{ll}
\displaystyle(\de-i\omega_{n})(\de+i\omega_{n})(\de^{2}-\omega_{n}^{2})\widetilde{w}_{1,n}=F_{1,n}&\text{ in }(0,\alpha)\cup(\beta,l)
\end{array}
\end{equation}
where  we denoted by $\displaystyle\de=\frac{\ud}{\ud x}$.

On the interval $(0,\alpha)$, by solving the first linear equation of~\eqref{pwkv44} we get
\begin{equation*}
(\de+i\omega_{n})(\de^{2}-\omega_{n}^{2})\widetilde{w}_{1,n}=C_{1}'\e^{i\omega_{n}(x-\alpha)}+\int_{\alpha}^{x}\e^{i\omega_{n}(x-s)}F_{1,n}(s)\,\ud s,
\end{equation*}
and
\begin{equation}\label{pwkv46}
(\de^{2}-\omega_{n}^{2})\widetilde{w}_{1,n}=C_{2}'\e^{-i\omega_{n}(x-\alpha)}+\frac{C_{1}'}{\omega_{n}}\sin(\omega_{n}(x-\alpha))+\frac{1}{\omega_{n}}\int_{\alpha}^{x}\sin(\omega_{n}(x-s))F_{1,n}(s)\,\ud s,
\end{equation}
where
\begin{equation}\label{pwkv51}
C_{1}'=\widetilde{w}_{1,n}'''(\alpha^{-})+i\omega_{n} \widetilde{w}_{1,n}''(\alpha^{-})-\omega_{n}^{2}\widetilde{w}_{1,n}'(\alpha)-i\omega_{n}^{3}\widetilde{w}_{1,n}(\alpha),
\end{equation}
and
\begin{equation}\label{pwkv52}
C_{2}'=\widetilde{w}_{1,n}''(\alpha^{-})-\omega_{n}^{2}\widetilde{w}_{1,n}(\alpha).
\end{equation}
We further solve~\eqref{pwkv46} and using the boundary conditions $\widetilde{w}_{1,n}(0)=\widetilde{w}_{1,n}'(0)=0$ to get
\begin{equation}\label{pwkv48}
\begin{split}
(\de-\omega_{n})\widetilde{w}_{1,n}=\frac{C_{2}'}{(1-i)\omega_{n}}\left[\e^{-i\omega_{n}(x-\alpha)}-\e^{-\omega_{n}(x-i\alpha)}\right]
\\
+\frac{C_{1}'}{2\omega_{n}^{2}}\left[\sin(\omega_{n}(x-\alpha))-\cos(\omega_{n}(x-\alpha))+\e^{-\omega_{n} x}(\sin(\omega_{n}\alpha)+\cos(\omega_{n}\alpha))\right]
\\
+\frac{1}{\omega_{n}}\int_{0}^{x}\int_{\alpha}^{\tau}\e^{-\omega_{n}(x-\tau)}\sin(\omega_{n}(\tau-s))F_{1,n}(s)\,\ud s\,\ud\tau.
\end{split}
\end{equation}
Multiplying~\eqref{pwkv48} by $2\omega_{n}$ and taking $x=\alpha$, we have
\begin{equation}\label{pwkv53}
\begin{split}
2\omega_{n} \widetilde{w}_{1,n}'(\alpha)-2\omega_{n}^{2}w_{1,n}(\alpha)=(1+i)C_{2}'(1-\e^{-\alpha\omega_{n}}\e^{i\alpha\omega_{n}})
\\
+\frac{C_{1}'}{\omega_{n}}(\e^{-\alpha\omega_{n}}(\sin(\alpha\omega_{n})+\cos(\alpha\omega_{n}))-1)
\\
+2\int_{0}^{\alpha}\int_{\alpha}^{\tau}\e^{-\omega_{n}(\alpha-\tau)}\sin(\omega_{n}(\tau-s))F_{1,n}(s)\,\ud s\,\ud\tau.
\end{split}
\end{equation}
We substitute~\eqref{pwkv51} and~\eqref{pwkv52} into~\eqref{pwkv53} and let $n\,\longrightarrow\,+\infty$, then~\eqref{pwkv23} and~\eqref{pwkv25} yields
$$\lim_{n\rightarrow+\infty}\widetilde{w}_{1,n}''(\alpha^{-})=-2\lim_{n\rightarrow+\infty}\int_{0}^{\alpha}\int_{\alpha}^{\tau}\e^{-\omega_{n}(\alpha-\tau)}\sin(\omega_{n}(\tau-s))F_{1,n}(s)\,\ud s\,\ud\tau.$$
We argue that the above limit is zero by the following estimates
$$\left|\int_{0}^{\alpha}\int_{\alpha}^{\tau}\e^{-\omega_{n}(\alpha-\tau)}\sin(\omega_{n}(\tau-s))g_{1,n}(s)\,\ud s\,\ud\tau\right|\leq\alpha^{\frac{3}{2}}\int_{0}^{\alpha}|g_{1,n}(s)|^{2}\ud s\,\longrightarrow\,0,$$
and
\begin{equation*}
\begin{split}
\left|\int_{0}^{\alpha}\int_{\alpha}^{\tau}\e^{-\omega_{n}(\alpha-\tau)}\sin(\omega_{n}(\tau-s))\lambda_{n}f_{1,n}(s)\,\ud s\,\ud\tau\right|=
\\
\left|\lambda_{n}\int_{0}^{\alpha}\left(\int_{0}^{s}\e^{-\omega_{n}(\alpha-\tau)}\sin(\omega_{n}(\tau-s))\,\ud\tau\right)f_{1,n}(s)\,\ud s\right|=
\\
\left|\frac{\lambda_{n}\e^{-\alpha\omega_{n}}}{2\omega_{n}}\int_{0}^{\alpha}(\cos(\omega_{n} s)+\sin(\omega_{n} s)-\e^{\omega_{n} s})f_{1,n}(s)\,\ud s\right|\leq
\\
\left(\alpha\omega_{n}\e^{-\alpha\omega_{n}}+\frac{1}{2}-\frac{1}{2}\e^{-\alpha\omega_{n}}\right)\max_{[0,\alpha]}|f_{1,n}(s)|\,\longrightarrow\,0,
\end{split}
\end{equation*}
where we have used the fact that $g_{n}\,\longrightarrow\,0$ in $H$, $f_{n}\,\longrightarrow\,0$ in $V_{1}\hookrightarrow\mathrm{C}^{1}([0,l])$, and $\omega_{n}\,\longrightarrow\,+\infty$. Thus we have proved the first identity of~\eqref{pwkv43}.

On the interval $(\beta,l)$, by solving the first linear equation~\eqref{pwkv44} we get
\begin{equation}\label{pwkv45}
(\de+i\omega_{n})(\de^{2}-\omega_{n}^{2})\widetilde{w}_{1,n}=C_{1}\e^{i\omega_{n}(x-\beta)}+\int_{\beta}^{x}\e^{i\omega_{n}(x-s)}F_{1,n}(s)\,\ud s,
\end{equation}
\begin{equation}\label{pwkv47}
(\de^{2}-\omega_{n}^{2})\widetilde{w}_{1,n}=C_{2}\e^{-i\omega_{n}(x-\beta)}+\frac{C_{1}}{\omega_{n}}\sin(\omega_{n}(x-\beta))+\frac{1}{\omega_{n}}\int_{\beta}^{x}\sin(\omega_{n}(x-s))F_{1,n}(s)\,\ud s,
\end{equation}
where
\begin{equation}\label{pwkv54}
C_{1}=\widetilde{w}_{1,n}'''(\beta^{+})+i\omega_{n} \widetilde{w}_{1,n}''(\beta^{+})-\omega_{n}^{2}\widetilde{w}_{1,n}'(\beta)-i\omega_{n}^{3}\widetilde{w}_{1,n}(\beta),
\end{equation}
\begin{equation}\label{pwkv92}
C_{2}=\widetilde{w}_{1,n}''(\beta^{+})-\omega_{n}^{2}\widetilde{w}_{1,n}(\beta),
\end{equation}
and
\begin{equation}\label{pwkv49}
\begin{split}
(\de+\omega_{n})\widetilde{w}_{1,n}=\frac{C_{2}}{(1+i)\omega_{n}}\left[\e^{\omega_{n}(x-\beta)}-\e^{-i\omega_{n}(x-\beta)}\right]
\\
+\frac{C_{1}}{2\omega_{n}^{2}}\left[1-\cos(\omega_{n}(x-\beta))-\sin(\omega_{n}(x-\beta))\right]+C_{3}\e^{\omega_{n} (x-\beta)}
\\
+\frac{1}{\omega_{n}}\int_{\beta}^{x}\int_{\beta}^{\tau}\e^{\omega_{n}(x-\tau)}\sin(\omega_{n}(\tau-s))F_{1,n}(s)\,\ud s\,\ud\tau,
\end{split}
\end{equation}
where
\begin{equation}\label{pwkv57}
C_{3}=\widetilde{w}_{1,n}'(\beta)+\omega_{n} \widetilde{w}_{1,n}(\beta).
\end{equation}
Multiplying the relation~\eqref{pwkv49} by $\e^{-\omega_{n}(x-\beta)}$ and taking $x=l$ then $w_{1,n}'(l)=0$ leads to
\begin{equation}\label{pwkv91}
\begin{split}
\e^{-\omega_{n}(l-\beta)}\widetilde{w}_{1,n}(l)=\frac{C_{2}(1-i)}{2\omega_{n}^{2}}\left[1-\e^{-\omega_{n}(1+i)(l-\beta)}\right]
\\
+\frac{C_{1}}{2\omega_{n}^{3}}\e^{-\omega_{n}(l-\beta)}\left[1-\cos(\omega_{n}(l-\beta))-\sin(\omega_{n}(l-\beta))\right]+\frac{C_{3}}{\omega_{n}}
\\
+\frac{1}{\omega_{n}^{2}}\int_{\beta}^{l}\int_{\beta}^{\tau}\e^{\omega_{n}(2\beta-l-\tau)}\sin(\omega_{n}(\tau-s))F_{1,n}(s)\,\ud s\,\ud\tau,
\end{split}
\end{equation}
where the last term satisfy
\begin{equation}\label{pwkv106}
\int_{\beta}^{l}\int_{\beta}^{\tau}\e^{\omega_{n}(2\beta-l-\tau)}\sin(\omega_{n}(\tau-s))F_{1,n}(s)\,\ud s\,\ud\tau\,\longrightarrow\,0.
\end{equation}
Indeed, we have
$$
\left|\int_{\beta}^{l}\int_{\beta}^{\tau}\e^{\omega_{n}(2\beta-l-\tau)}\sin(\omega_{n}(\tau-s))g_{1,n}(s)\,\ud s\,\ud\tau\right|\leq(l-\beta)^{\frac{3}{2}}\left(\int_{\beta	}^{l}|g_{1,n}(s)|^{2}\ud s\right)^{\frac{1}{2}}\longrightarrow\,0,
$$
and
$$
\left|\int_{\beta}^{l}\!\!\!\int_{\beta}^{\tau}\!\!\e^{\omega_{n}(2\beta-l-\tau)}\sin(\omega_{n}(\tau-s))\lambda_{n}f_{1,n}(s)\,\ud s\,\ud\tau\right|\leq C|\lambda_{n}|\e^{-\omega_{n}(l-\beta)}\max_{[\beta,l]}|f_{1,n}(s)|\longrightarrow\,0.
$$
Substitute the expression of $C_{1}$, $C_{2}$ and $C_{3}$ in~\eqref{pwkv54},~\eqref{pwkv92} and~\eqref{pwkv57} respectively into~\eqref{pwkv91} we obtain
\begin{equation*}
\begin{split}
\widetilde{w}_{1,n}''(\beta^{+})\left[1-\left(1+i(\cos(\omega_{n}(l-\beta))+\sin(\omega_{n}(l-\beta))-1)\right)\e^{-\omega_{n}(l-\beta)}\right]
\\
=\omega_{n}^{-1}\widetilde{w}_{1,n}'''(\beta^{+})\left[-1+\cos(\omega_{n}(l-\beta))+\sin(\omega_{n}(l-\beta))\right]\e^{-\omega_{n}(l-\beta)}
\\
+\omega_{n}\widetilde{w}_{1,n}'(\beta)\left[-(1+i)+\left(1-\cos(\omega_{n}(l-\beta))-\sin(\omega_{n}(l-\beta))\right)\e^{-\omega_{n}(l-\beta)}\right]
\\
+\omega_{n}^{2}\widetilde{w}_{1,n}(\beta)\left[-i+\left(i-(1+i)\cos(\omega_{n}(l-\beta))\right)\e^{-\omega(l-\beta)}\right]-(1+i)\omega_{n}^{2}\widetilde{w}_{1,n}(l)\e^{-\omega_{n}(l-\beta)}
\\
-(1+i)\int_{\beta}^{l}\int_{\beta}^{\tau}\e^{\omega_{n}(2\beta-l-\tau)}\sin(\omega_{n}(\tau-s))F_{1,n}(s)\,\ud s\,\ud\tau.
\end{split}
\end{equation*}
Since $|\lambda_{n}w_{1,n}(l)|$ is bounded by~\eqref{pwkv42} then the second statement of~\eqref{pwkv43} follows from~\eqref{pwkv23}, \eqref{pwkv25} in Lemma~\ref{pwkv111},~\eqref{pwkv106} and the fact that $\omega_{n}\,\longrightarrow\,+\infty$ as $n\,\longrightarrow\,+\infty$.

This leads form~\eqref{pwkv102} to
\begin{equation}\label{pwkv50}
|w_{1,n}''(l)|^{2}\,\longrightarrow\,0\quad\text{and}\quad|\lambda_{n}w_{2,n}(l)|^{2}\,\longrightarrow\,0
\end{equation}
and form~\eqref{pwkv103} we get
$$
\int_{l}^{L}|\lambda_{n}.w_{2,n}|^{2}+3\int_{l}^{L}|w_{2,n}'|^{2}\,\ud x\longrightarrow\,2
$$
which combined with~\eqref{pwkv31} imply that
$$
\int_{l}^{L}|\lambda_{n}.w_{2,n}|^{2}\,\longrightarrow\,\frac{1}{2}\quad\text{and}\quad\int_{l}^{L}|w_{2,n}'|^{2}\,\ud x\longrightarrow\,\frac{1}{2}.
$$
Then from~\eqref{pwkv42} and~\eqref{pwkv43} we obtain
\begin{equation}\label{pwkv55}
(L-l)|w_{2,n}'(l)|^{2}\,\longrightarrow\,1.
\end{equation}

In what follows, and in order to achieve our proof we will try to obtain a contradiction with~\eqref{pwkv55}. 
\\
We take $x=l$ in the relation~\eqref{pwkv47} then we find
\begin{equation}\label{pwkv97}
\widetilde{w}_{1,n}''(l)=\omega_{n}^{2}\widetilde{w}_{1,n}(l)+C_{2}\e^{-i\omega_{n}(l-\beta)}+\frac{C_{1}}{\omega_{n}}\sin(\omega_{n}(l-\beta))+\frac{1}{\omega_{n}}\int_{\beta}^{l}\sin(\omega_{n}(l-s))F_{1,n}(s)\,\ud s,
\end{equation}
where the last term verify
\begin{equation}\label{pwkv105}
\frac{1}{\omega_{n}}\left|\int_{\beta}^{l}\sin(\omega_{n}(l-s))F_{1,n}(s)\,\ud s\right|\,\longrightarrow\,0.
\end{equation}
Indeed, we have
$$
\left|\int_{\beta}^{l}\sin(\omega_{n}(l-s))g_{1,n}(s)\,\ud s\right|\leq(l-\beta)^{\frac{1}{2}}\left(\int_{\beta	}^{l}|g_{1,n}(s)|^{2}\ud s\right)^{\frac{1}{2}}\longrightarrow\,0,
$$
and
$$
\left|\int_{\beta}^{l}\sin(\omega_{n}(l-s))f_{1,n}(s)\,\ud s\right|\leq C\max_{[\beta,l]}|f_{1,n}(s)|\longrightarrow\,0.
$$
Taking again $x=l$ in the relation~\eqref{pwkv45},  and using~\eqref{pwkv97}, $w_{1,n}'(l)=0$ and the transmission condition $w_{1,n}'''(l)+w_{2,n}'(l)=0$ then we obtain
\begin{equation}\label{pwkv98}
\begin{split}
\widetilde{w}_{1,n}'''(l)=-\omega_{n}w_{2,n}'(l)=-i\omega_{n} C_{2}\e^{-i\omega_{n}(l-\beta)}+C_{1}\cos(\omega_{n}(l-\beta))
\\
-i\int_{\beta}^{l}\sin(\omega_{n}(l-s))F_{1,n}(s)\,\ud s+\int_{\beta}^{l}\e^{i\omega_{n}(l-s)}F_{1,n}(s)\,\ud s,
\end{split}
\end{equation}
where the last term verify
\begin{equation}\label{pwkv104}
\frac{1}{\omega_{n}}\left|\int_{\beta}^{l}\e^{i\omega_{n}(l-s)}F_{1,n}(s)\,\ud s\right|\,\longrightarrow\,0.
\end{equation}
Indeed, we have
$$
\left|\int_{\beta}^{l}\e^{i\omega_{n}(l-s)}g_{1,n}(s)\,\ud s\right|\leq(l-\beta)^{\frac{1}{2}}\left(\int_{\beta	}^{l}|g_{1,n}(s)|^{2}\ud s\right)^{\frac{1}{2}}\longrightarrow\,0,
$$
and
$$
\left|\int_{\beta}^{l}\e^{i\omega_{n}(l-s)}f_{1,n}(s)\,\ud s\right|\leq C\max_{[\beta,l]}|f_{1,n}(s)|\longrightarrow\,0.
$$
We substitute the expression of $C_{1}$ and $C_{2}$ as defined in ~\eqref{pwkv54} and~\eqref{pwkv92} respectively into the relation~\eqref{pwkv98} then we find
\begin{equation}\label{pwkv56}
\begin{split}
w_{2,n}'(l)=-(w_{1,n}'''(\beta^{+})+\omega_{n}^{2}w_{1,n}(\beta))\cos(\omega_{n}(l-\beta))+(\omega_{n}w_{1,n}''(\beta^{+})-\omega_{n}^{3}w_{1,n}(\beta))\sin(\omega_{n}(l-\beta))
\\
+\frac{i}{\omega_{n}}\int_{\beta}^{l}\sin(\omega_{n}(l-s))F_{1,n}(s)\,\ud s-\frac{1}{\omega_{n}}\int_{\beta}^{l}\e^{i\omega_{n}(l-s)}F_{1,n}(s)\,\ud s.
\end{split}
\end{equation}
Finally, since the terms in right hand side of~\eqref{pwkv56} tends to zero as $n\,\longrightarrow\,+\infty$ by using~\eqref{pwkv23}, \eqref{pwkv25} in Lemma~\ref{pwkv111},~\eqref{pwkv43},~\eqref{pwkv105} and~\eqref{pwkv104} then we find that $w_{2,n}'(l)\,\longrightarrow\,0$. Hence we proved the promised contradiction. And this complete the proof.
\section{Damping arising from the longitudinal motion}\label{pwkv95}
The purpose of this section is to prove the first part of Theorem~\ref{pwkv121}. We need only to verify the condition for a semigroup of contractions on a Hilbert space being polynomially stable (see~\cite{BT}), i.e.,
\begin{equation}\label{pwkv61}
\sup\{|\lambda|^{-1}.\|(i\lambda-\mathcal{A}_{2})^{-1}\|;\quad\lambda\in\R\}<+\infty.
\end{equation}
We will argue by contradiction, thus we suppose that~\eqref{pwkv61} is not true. By the continuity of the resolvent and the resonance theorem, there exist $\lambda_{n}\in\R$, $((u_{1,n},u_{2,n}),(v_{1,n},v_{2,n}))\in\mathcal{D}(\mathcal{A}_{2})$, $n=1,2,\ldots$, such that
\begin{equation}\label{pwkv62}
\|((u_{1,n},u_{2,n}),(v_{1,n},v_{2,n}))\|_{\mathcal{H}_{2}}=1,\qquad\lambda_{n}\,\longrightarrow\,+\infty
\end{equation}
and 
\begin{equation}\label{pwkv63}
\lambda_{n}(i\lambda_{n}-\mathcal{A}_{2})((u_{1,n},u_{2,n}),(v_{1,n},v_{2,n}))=((f_{1,n},f_{2,n}),(g_{1,n},g_{2,n}))\,\longrightarrow\,0\;\text{in}\;\mathcal{H}_{2},
\end{equation}
which mean
\begin{eqnarray}
\lambda_{n}.(i\lambda_{n}u_{1,n}-v_{1,n})=f_{1,n}\,\longrightarrow\,0\;\text{in}\; H^{1}(0,l)\label{pwkv64}
\\
\lambda_{n}.(i\lambda_{n}u_{2,n}-v_{2,n})=f_{2,n}\,\longrightarrow\,0\;\text{in}\; H^{1}(l,L)\label{pwkv65}
\\
\lambda_{n}.(i\lambda_{n}v_{1,n}-T_{n}')=g_{1,n}\,\longrightarrow\,0\;\text{in}\; L^{2}(0,l)\label{pwkv66}
\\
\lambda_{n}.(i\lambda_{n}v_{2,n}+u_{2,n}'''')=g_{2,n}\,\longrightarrow\,0\;\text{in}\; L^{2}(l,L)\label{pwkv67}
\end{eqnarray}
where
$$
T_{n}=u_{1,n}'+D_{a}v_{1,n}'.
$$

We first consider~\eqref{pwkv64} and~\eqref{pwkv66} on the interval $(\alpha,\beta)$. From~\eqref{pwkv63}, we obtain
\begin{equation}\label{pwkv68}
\lambda_{n}\int_{\alpha}^{\beta}\!\!a|v_{1,n}'|^{2}\ud x=\re\left\langle\lambda_{n}(i\lambda_{n}-\mathcal{A}_{2})(u_{1,n},u_{2,n},v_{1,n},v_{2,n}),(u_{1,n},u_{2,n},v_{1,n},v_{2,n}) \right\rangle_{\mathcal{H}_{2}}\longrightarrow0
\end{equation}
which imply that
\begin{equation}\label{pwkv69}
\|\lambda_{n}^{\frac{1}{2}}v_{1,n}'\|_{L^{2}(\alpha,\beta)}\,\longrightarrow\,0,\quad\text{and}\quad \|\lambda_{n}^{\frac{3}{2}}u_{1,n}'\|_{L^{2}(\alpha,\beta)}\,\longrightarrow\,0.
\end{equation}
Thus we also have
\begin{equation}\label{pwkv70}
\|\lambda_{n}^{\frac{1}{2}}T_{n}\|_{L^{2}(\alpha,\beta)}\,\longrightarrow\,0.
\end{equation}
The rest of the proof depend on the following lemma.
\begin{lem}
The functions $u_{1,n}$ and $T_{n}$ have the following properties:
\begin{eqnarray}
&|\lambda_{n}^{\frac{1}{2}}T_{n}(\alpha^{+})|\,\longrightarrow\,0,&\qquad |\lambda_{n}^{\frac{1}{2}}T_{n}(\beta^{-})|\,\longrightarrow\,0,\label{pwkv71}
\\
&|\lambda_{n}^{\frac{3}{2}}u_{1,n}(\alpha^{+})|\,\longrightarrow\,0,&\qquad |\lambda_{n}^{\frac{3}{2}}u_{1,n}(\beta^{-})|\,\longrightarrow\,0.\label{pwkv72}
\end{eqnarray}
\end{lem}
\begin{pr}
From~\eqref{pwkv64} we have
\begin{equation}\label{pwkv73}
\lambda_{n}^{-1}.\|\psi v_{1,n}\|_{H^{1}(\alpha,\beta)}=O(1),\qquad\forall\,\psi\in\mathrm{C}^{\infty}(0,l),\;\supp(\psi)\subset(\alpha,\beta).
\end{equation}
Equations~\eqref{pwkv66},~\eqref{pwkv70} and~\eqref{pwkv73} imply that
\begin{equation}\label{pwkv74}
\lambda_{n}\int_{\alpha}^{\beta}\psi|v_{1,n}|^{2}\,\ud x\,\longrightarrow\,0,\qquad\forall\,\psi\in\mathrm{C}^{\infty}(0,l),\;\supp(\psi)\subset(\alpha,\beta).
\end{equation}
Applying the interpolation theorem involving subdomains~\cite[Theorem 4.23]{A} we find that~\eqref{pwkv69} and~\eqref{pwkv74} imply
\begin{equation}\label{pwkv75}
\lambda_{n}^{\frac{1}{2}}v_{1,n}\,\longrightarrow\,0\qquad\text{in }H^{1}(\alpha,\beta).
\end{equation}
We take the inner product of~\eqref{pwkv66} with $v_{1,n}$ in $L^{2}(\alpha,\beta)$ to obtain
\begin{equation}\label{pwkv76}
\begin{split}
\int_{\alpha}^{\beta}g_{1,n}.\overline{v}_{1,n}\,\ud x=i\lambda_{n}^{2}\int_{\alpha}^{\beta}|v_{1,n}|^{2}\,\ud x+\int_{\alpha}^{\beta}\lambda_{n}^{\frac{1}{2}}T_{n}.\lambda_{n}^{\frac{1}{2}}\overline{v}_{1,n}'\,\ud x
\\
+\lambda_{n}T_{n}(\alpha^{+})\overline{v}_{1,n}(\alpha^{+})-\lambda_{n}T_{n}(\beta^{-})\overline{v}_{1,n}(\beta^{-}).
\end{split}
\end{equation}
Using~\eqref{pwkv66},~\eqref{pwkv70} and~\eqref{pwkv75} we obtain for the third term in the right hand side of~\eqref{pwkv76} that
\begin{eqnarray}\label{pwkv77}
\lambda_{n}.|T_{n}(\alpha^{+})|.|v_{1,n}(\alpha^{+})|\!\!\!\!\!&\leq&\!\!\!C\lambda_{n}\|T_{n}\|_{L^{2}(\alpha,\beta)}^{\frac{1}{2}}.\|T_{n}'\|_{L^{2}(\alpha,\beta)}^{\frac{1}{2}}.\|v_{1,n}\|_{L^{2}(\alpha,\beta)}^{\frac{1}{2}}.\|v_{1,n}'\|_{L^{2}(\alpha,\beta)}^{\frac{1}{2}}\nonumber
\\
&=&\!\!\!C\|\lambda_{n}^{\frac{1}{2}}T_{n}\|_{L^{2}(\alpha,\beta)}^{\frac{1}{2}}.\|\lambda_{n}v_{1,n}\|_{L^{2}(\alpha,\beta)}.\|\lambda_{n}^{\frac{1}{2}}v_{1,n}'\|_{L^{2}(\alpha,\beta)}^{\frac{1}{2}}+o(1)\nonumber
\\
&=& o(1)(1+\|\lambda_{n}v_{1,n}\|_{L^{2}(\alpha,\beta)}).
\end{eqnarray}
And similarly we can show that
\begin{equation}\label{pwkv78}
\lambda_{n}.|T_{n}(\beta^{-})|.|v_{1,n}(\beta^{-})|\leq o(1)(1+\|\lambda_{n}v_{1,n}\|_{L^{2}(\alpha,\beta)}).
\end{equation}
Now~\eqref{pwkv76},~\eqref{pwkv77} and~\eqref{pwkv78} leads to
\begin{equation}\label{pwkv79}
\|\lambda_{n}v_{1,n}\|_{L^{2}(\alpha,\beta)}\,\longrightarrow\,0.
\end{equation}
Next we multiply~\eqref{pwkv66} by $\lambda_{n}^{-\frac{1}{2}}(\beta-x)\overline{T}_{n}$ and take inner product in $L^{2}(\alpha,\beta)$ to get
\begin{equation}\label{pwkv80}
\re\left\langle i\lambda_{n}^{\frac{3}{2}} v_{1,n},(\beta-x)T_{n}\right\rangle_{L^{2}(\alpha,\beta)}-\frac{\lambda_{n}^{\frac{1}{2}}}{2}(\beta-\alpha)|T_{n}(\alpha^{+})|^{2}-\frac{\lambda_{n}^{\frac{1}{2}}}{2}\|T_{n}\|_{L^{2}(\alpha,\beta)}^{2}\,\longrightarrow\,0,
\end{equation}
and by multiplying ~\eqref{pwkv66} by $\lambda_{n}^{-\frac{1}{2}}(\alpha-x)\overline{T}_{n}$ and take inner product in $L^{2}(\alpha,\beta)$ we get
\begin{equation}\label{pwkv81}
\re\left\langle i\lambda_{n}^{\frac{3}{2}}v_{1,n},(\alpha-x)T_{n}\right\rangle_{L^{2}(\alpha,\beta)}-\frac{\lambda_{n}^{\frac{1}{2}}}{2}(\beta-\alpha)|T_{n}(\beta^{-})|^{2}-\frac{\lambda_{n}^{\frac{1}{2}}}{2}\|T_{n}\|_{L^{2}(\alpha,\beta)}^{2}\,\longrightarrow\,0.
\end{equation}
Since the first and the third terms of~\eqref{pwkv80} and~\eqref{pwkv81} converge to zero by~\eqref{pwkv70} and~\eqref{pwkv79} then~\eqref{pwkv71} yields. On the other hand by~\eqref{pwkv64} and~\eqref{pwkv75} we obtain
\begin{equation}\label{pwkv82}
\lambda_{n}^{\frac{3}{2}}.\|u_{n}\|_{L^{2}(\alpha,\beta)}\,\longrightarrow\,0,
\end{equation}
then by~\eqref{pwkv69} and~\eqref{pwkv82} we get
\begin{equation}\label{pwkv60}
\lambda_{n}^{\frac{3}{2}}.\|u_{n}\|_{H^{1}(\alpha,\beta)}\,\longrightarrow\,0,
\end{equation}
hence~\eqref{pwkv72} hold from the Sobolev embedding inequalities.
\end{pr}

Using now the continuity conditions at $x=\alpha$ and $x=\beta$, we arrive from~\eqref{pwkv71} and~\eqref{pwkv72} at
\begin{equation}\label{pwkv83}
\begin{split}
|\lambda_{n}^{\frac{1}{2}}u_{1,n}'(\alpha^{-})|\,\longrightarrow\,0,&\qquad |\lambda_{n}^{\frac{1}{2}}u_{1,n}'(\beta^{+})|\,\longrightarrow\,0,
\\
|\lambda_{n}^{\frac{3}{2}}u_{1,n}(\alpha^{-})|\,\longrightarrow\,0,&\qquad |\lambda_{n}^{\frac{3}{2}}u_{1,n}(\beta^{+})|\,\longrightarrow\,0.
\end{split}
\end{equation}
We consider now~\eqref{pwkv64}-\eqref{pwkv67} on the intervals $(0,\alpha)$, $(\beta,l)$ and $(l,L)$, then by replacing~\eqref{pwkv64} and~\eqref{pwkv65} respectively into~\eqref{pwkv66} and~\eqref{pwkv67} we obtain
\begin{eqnarray}
-\lambda_{n}^{2}u_{1,n}-u_{1,n}''=\lambda_{n}^{-1}g_{1,n}+if_{1,n}&\text{in}&(0,\alpha)\cup(\beta,l),\label{pwkv84}
\\
-\lambda_{n}^{2}u_{2,n}+u_{2,n}''''=\lambda_{n}^{-1}g_{2,n}+if_{2,n}&\text{in}&(l,L).\label{pwkv85}
\end{eqnarray}
Take the inner product of~\eqref{pwkv84} with $\displaystyle\frac{x}{2}u_{1,n}'$ in $L^{2}(0,\alpha)$ and with $\displaystyle\frac{(x-L)}{2}u_{1,n}'$ in $L^{2}(\beta,l)$ and the inner product of~\eqref{pwkv85} with $(L-x)u_{2,n}'$ in $L^{2}(l,L)$. A straight forward calculation shows that the real part of this inner products leads to the following
\begin{equation}\label{pwkv86}
\frac{1}{2}\int_{0}^{\alpha}|\lambda_{n}u_{1,n}|^{2}\,\ud x+\frac{1}{2}\int_{0}^{\alpha}|u_{1,n}'|^{2}\,\ud x-\frac{\alpha}{2}|\lambda_{n}u_{1,n}(\alpha)|^{2}-\frac{\alpha}{2}|u_{1,n}'(\alpha^{-})|^{2}\,\longrightarrow\,0,
\end{equation}
\begin{equation}\label{pwkv87}
\begin{split}
\frac{1}{2}\int_{\beta}^{l}|\lambda_{n}u_{1,n}|^{2}\,\ud x+\frac{1}{2}\int_{\beta}^{l}|u_{1,n}'|^{2}\,\ud x+\frac{(L-l)}{2}\left(|\lambda_{n}u_{1,n}(l)|^{2}+|u_{1,n}'(l)|^{2}\right)
\\
-\frac{(L-\beta)}{2}\left(|\lambda_{n}u_{1,n}(\beta)|^{2}+|u_{1,n}'(\beta^{+})|^{2}\right)\,\longrightarrow\,0,
\end{split}
\end{equation}
and
\begin{equation}\label{pwkv88}
\frac{1}{2}\int_{l}^{L}|\lambda_{n}u_{2,n}|^{2}\,\ud x+\frac{3}{2}\int_{l}^{L}|u_{2,n}''|^{2}\,\ud x-\frac{(L-l)}{2}\left(|\lambda_{n}u_{2,n}(l)|^{2}+|u_{2,n}''(l)|^{2}\right)\,\longrightarrow\,0.
\end{equation}
Similar as the previous section we prove also that
\begin{equation}\label{pwkv59}
\begin{split}
\lim_{n\rightarrow+\infty}\int_{0}^{\alpha}|\lambda_{n}.u_{1,n}|^{2}\,\ud x+\int_{\beta}^{l}|\lambda_{n}.u_{1,n}|^{2}\,\ud x+\int_{l}^{L}|\lambda_{n}.u_{2,n}|^{2}\,\ud x=
\\
\lim_{n\rightarrow+\infty}\int_{0}^{\alpha}|u_{1,n}'|^{2}\,\ud x+\int_{\beta}^{l}|u_{1,n}'|^{2}\,\ud x+\int_{l}^{L}|u_{2,n}''|^{2}\,\ud x=\frac{1}{2}
\end{split}
\end{equation}
where~\eqref{pwkv60} is taken into account.
\\
Summing now~\eqref{pwkv86}-\eqref{pwkv88} and using~\eqref{pwkv83} and~\eqref{pwkv59} then we find
\begin{equation}\label{pwkv58}
\begin{split}
(L-l)\left(|u_{2,n}''(l)|^{2}-|u_{1,n}'(l)|^{2}\right)+2\int_{\beta}^{l}|u_{1,n}'|^{2}\,\ud x\,\longrightarrow\,2.
\end{split}
\end{equation}

Equation~\eqref{pwkv84} is rewritten now in $(\beta,l)$ as follows
\begin{equation}\label{pwkv89}
(\de-i\lambda_{n})(\de+i\lambda_{n})u_{1,n}=F_{1,n}\quad\text{in  }(\beta,l),
\end{equation}
where we have denoted by $\displaystyle\de=\frac{\ud}{\ud x}$ and $F_{1,n}=-\lambda_{n}^{-1}g_{1,n}-if_{1,n}$.
\\
Solving the first order equation of~\eqref{pwkv89}, we have
\begin{equation}\label{pwkv90}
(\de+i\lambda_{n})u_{1,n}(x)=K_{1}'\e^{i\lambda_{n}(x-\beta)}+\int_{\beta}^{x}\e^{i\lambda_{n}(x-s)}F_{1,n}(s)\,\ud s,
\end{equation}
where
\begin{equation}\label{pwkv93}
K_{1}'=u_{1,n}'(\beta^{+})+i\lambda_{n}u_{1,n}(\beta).
\end{equation}
Resolving the first order equation of~\eqref{pwkv90}, we have
\begin{equation}\label{pwkv109}
u_{1,n}(x)=\frac{K_{1}'}{\lambda_{n}}\sin(\lambda_{n}(x-\beta))+K_{2}'\e^{-i\lambda_{n}(x-\beta)}+\frac{1}{\lambda_{n}}\int_{\beta}^{x}\sin(\lambda_{n}(x-s))F_{1,n}(s)\,\ud s,
\end{equation}
where
\begin{equation}\label{pwkv110}
K_{2}'=u_{1,n}(\beta).
\end{equation}
Multiplying~\eqref{pwkv109} by $\lambda_{n}$ and taking $x=l$, we find
\begin{equation}\label{pwkv112}
\lambda_{n}u_{1,n}(l)=K_{1}'\sin(\lambda_{n}(l-\beta))+\lambda_{n}K_{2}'\e^{-i\lambda_{n}(l-\beta)}+\int_{\beta}^{l}\sin(\lambda_{n}(l-s))F_{1,n}(s)\,\ud s.
\end{equation}
Since $K_{1}'$ and $\lambda_{n}K_{2}'$, as defined in~\eqref{pwkv93} and~\eqref{pwkv110} respectively, converge to zero using~\eqref{pwkv83}, and with integrating by parts we have
$$
\left|\int_{\beta}^{l}\sin(\lambda_{n}(l-s))g_{1,n}(s)\,\ud s\right|\leq\sqrt{(l-\beta)}\left(\int_{\beta}^{l}|g_{1,n}(s)|^{2}\,\ud s\right)^{\frac{1}{2}}\,\longrightarrow\,0
$$
and
$$
\left|\int_{\beta}^{l}\sin(\lambda_{n}(l-s))f_{1,n}(s)\,\ud s\right|\leq\sqrt{(l-\beta)}\left(\int_{\beta}^{l}|f_{1,n}(s)|^{2}\,\ud s\right)^{\frac{1}{2}}\,\longrightarrow\,0,
$$
then the right hand side of~\eqref{pwkv112} converge to zero. This leads to
\begin{equation}\label{pwkv113}
|\lambda_{n}u_{1,n}(l)|\,\longrightarrow\,0.
\end{equation}
We take the derivative of~\eqref{pwkv109} at the point $x=l$, we obtain
\begin{equation}\label{pwkv114}
u_{1,n}'(l)=K_{1}'\cos(\lambda_{n}(l-\beta))-i\lambda_{n}K_{2}'\e^{-i\lambda_{n}(l-\beta)}+\int_{\beta}^{l}\cos(\lambda_{n}(l-s))F_{1,n}(s)\,\ud s.
\end{equation}
We use the same arguments as previously to prove that
\begin{equation}\label{pwkv115}
|u_{1,n}'(l)|\,\longrightarrow\,0.
\end{equation}

Equation~\eqref{pwkv85} is rewritten now as follows
\begin{equation}\label{pwkv108}
(\de-i\omega_{n})(\de+i\omega_{n})(\de^{2}-\omega_{n}^{2})u_{2,n}=F_{2,n}\quad\text{in  }(l,L)
\end{equation}
where we denoted by $F_{2,n}=\lambda_{n}^{-1}g_{2,n}+if_{2,n}$ and $\displaystyle\omega_{n}=\sqrt{\lambda_{n}}$.
\\
Solving the first linear equation of~\eqref{pwkv108} we get
\begin{equation*}
(\de+i\omega_{n})(\de^{2}-\omega_{n}^{2})u_{2,n}=K_{1}\e^{i\omega_{n}(x-l)}+\int_{l}^{x}\e^{i\omega_{n}(x-s)}F_{2,n}(s)\,\ud s,
\end{equation*}
\begin{equation*}
(\de^{2}-\omega_{n}^{2})u_{2,n}=K_{2}\e^{-i\omega_{n}(x-l)}+\frac{K_{1}}{\omega_{n}}\sin(\omega_{n}(x-l))+\frac{1}{\omega_{n}}\int_{l}^{x}\sin(\omega_{n}(x-s))F_{2,n}(s)\,\ud s,
\end{equation*}
where
\begin{eqnarray}\label{pwkv100}
K_{1}=u_{2,n}'''(l)+i\omega_{n} u_{2,n}''(l)-i\omega_{n}^{3}u_{2,n}(l),
\end{eqnarray}
\begin{equation}\label{pwkv101}
K_{2}=u_{2,n}''(l)-\omega_{n}^{2}u_{2,n}(l)
\end{equation}
and using now the boundary conditions $u_{2,n}(L)=u_{2,n}'(L)=0$ then we get
\begin{equation}\label{pwkv107}
\begin{split}
(\de+\omega_{n})u_{2,n}=\frac{K_{2}}{(1+i)\omega_{n}}\left[\e^{-\omega_{n}(L-x+i(L-l))}-\e^{-i\omega_{n}(x-l)}\right]
\\
+\frac{K_{1}}{2\omega_{n}^{2}}\left[-\cos(\omega_{n}(x-l))-\sin(\omega_{n}(x-l))+\e^{\omega_{n} (x-L)}(\cos(\omega_{n}(L-l))+\sin(\omega_{n}(L-l)))\right]
\\
+\frac{1}{\omega_{n}}\int_{L}^{x}\int_{l}^{\tau}\e^{\omega_{n}(x-\tau)}\sin(\omega_{n}(\tau-s))F_{2,n}(s)\,\ud s\,\ud\tau.
\end{split}
\end{equation}
Multiplying ~\eqref{pwkv107} by $-2\omega_{n}$ and take $x=l$, we arrive at
\begin{equation}\label{pwkv116}
\begin{split}
-2\omega_{n}^{2}u_{2,n}(l)=(1-i)K_{2}\left[1-\e^{-\omega_{n}(1+i)(L-l)}\right]+2\int_{l}^{L}\!\!\!\int_{l}^{\tau}\e^{\omega_{n}(l-\tau)}\sin(\omega_{n}(\tau-s))F_{2,n}(s)\,\ud s\,\ud\tau
\\
+\frac{K_{1}}{\omega_{n}}\left[1-\e^{-\omega_{n} (L-l)}(\cos(\omega_{n}(L-l))+\sin(\omega_{n}(L-l)))\right].
\end{split}
\end{equation}
We substitute~\eqref{pwkv100} and~\eqref{pwkv101} into~\eqref{pwkv116} and let $n\,\longrightarrow\,+\infty$, then~\eqref{pwkv113} and~\eqref{pwkv115} yields
$$\lim_{n\rightarrow+\infty}u_{2,n}''(l)=-2\lim_{n\rightarrow+\infty}\int_{l}^{L}\int_{l}^{\tau}\e^{\omega_{n}(l-\tau)}\sin(\omega_{n}(\tau-s))F_{2,n}(s)\,\ud s\,\ud\tau.$$
We argue that the above limit is zero by the following estimates
$$
\left|\int_{l}^{L}\!\!\!\int_{l}^{\tau}\e^{\omega_{n}(l-\tau)}\sin(\omega_{n}(\tau-s))g_{2,n}(s)\,\ud s\,\ud\tau\right|\leq(L-l)^{\frac{3}{2}}\left(\int_{l}^{L}|g_{2,n}(s)|^{2}\ud s\right)^{\frac{1}{2}}\,\longrightarrow\,0,
$$
and
\begin{equation*}
\begin{split}
\left|\int_{l}^{L}\!\!\!\int_{l}^{\tau}\e^{\omega_{n}(l-\tau)}\sin(\omega_{n}(\tau-s))f_{2,n}(s)\,\ud s\,\ud\tau\right|\leq(L-l)^{\frac{3}{2}}\left(\int_{l}^{L}|f_{2,n}(s)|^{2}\ud s\right)^{\frac{1}{2}}\,\longrightarrow\,0.
\end{split}
\end{equation*}
Hence we get that
\begin{equation}\label{pwkv117}
|u_{2,n}''(l)|\,\longrightarrow\,0.
\end{equation}
Returning now to~\eqref{pwkv58} and  using~\eqref{pwkv115} and~\eqref{pwkv117} then we obtain that
\begin{equation}\label{pwkv118}
\int_{\beta}^{l}|u_{1,n}'(x)|\,\ud x\,\longrightarrow\,1.
\end{equation}
Finally, we combine~\eqref{pwkv59} and~\eqref{pwkv118} to obtain a contradiction and this concludes the proof.
\section{Non exponential stability}\label{pwkv119}
The purpose of this section is to prove that the solutions of systems~\eqref{pwkv1} and~\eqref{pwkv2} are not exponentially stable. More precisely we will show that the resolvent $(i\lambda-\mathcal{A}_{1})$ and $(i\lambda-\mathcal{A}_{2})$ are not uniformly bounded with respect to $\lambda\in\R$.
\subsection{Case when the damping arising from the wave equation}
Let $\displaystyle\lambda=\lambda_{n}=\frac{4n^{2}\pi^{2}}{(L-l)^{2}}$, $n=1,2,\dots,$. We take $a$ constant in $(\alpha,\beta)$ and
\begin{equation*}
f(x)=f_{n}(x)=\left\{\begin{array}{cl}
\displaystyle\frac{1}{\lambda}\sin(\lambda(x-\alpha))&0<x<\alpha
\\
0&\alpha<x<L
\end{array}\right.
\end{equation*}
and
\begin{equation*}
g(x)=g_{n}(x)=\left\{\begin{array}{cl}
\displaystyle\cos(\lambda(x-\alpha))&0<x<\alpha
\\
0&\alpha<x<L.
\end{array}\right.
\end{equation*}
We set $f_{1,n},\,f_{2,n}$ the restriction of $f_{n}$ over the intervals $(0,l)$ and $(l,L)$ respectively. Similarly $g_{1,n},\,g_{2,n}$ are the restriction of $g_{n}$ over the intervals $(0,l)$ and $(l,L)$ respectively. In the intervals $(0,\alpha),\,(\alpha,\beta),\,(\beta,l)$ and $(l,L)$ we solve the resolvent equation
$$
(i\lambda_{n}-\mathcal{A}_{2})\left(\begin{array}{l}
u_{1,n}
\\
u_{2,n}
\\
v_{1,n}
\\
v_{2,n}
\end{array}\right)=\left(\begin{array}{l}
f_{1,n}
\\
f_{2,n}
\\
g_{1,n}
\\
g_{2,n}
\end{array}\right)
$$
where $(u,v)=(u_{n},v_{n})=((u_{1,n},u_{2,n}),(v_{1,n},v_{2,n}))\in\mathcal{D}(\mathcal{A}_{2})$.

For $x\in(0,\alpha)$, we have
\begin{equation}\label{pwkv123}
\left\{\begin{array}{ll}
i\lambda u-v=f&\text{in }(0,\alpha)
\\
i\lambda v-u''=g&\text{in }(0,\alpha)
\\
u(0)=0.&
\end{array}\right.
\end{equation}
Let $\displaystyle z_{\pm}=\frac{1}{2}(v(x)\pm u'(x))$ for all $x\in(0,\alpha)$. Then~\eqref{pwkv123} can be transformed into the first-order, diagonal and non homogeneous system in $(0,\alpha)$ as follow
$$
\left(\begin{array}{l}
z_{+}
\\
z_{-}
\end{array}\right)'=\left(\begin{array}{lr}
i\lambda&0
\\
0&-i\lambda
\end{array}\right)\left(\begin{array}{l}
z_{+}
\\
z_{-}
\end{array}\right)+\left(\begin{array}{r}
-1
\\
0
\end{array}\right)g:=A\left(\begin{array}{l}
z_{+}
\\
z_{-}
\end{array}\right)+\left(\begin{array}{r}
-1
\\
0
\end{array}\right)g.
$$
Using the boundary condition $v(0)=0=z_{+}(0)+z_{-}(0)$ we obtain the solution
$$
\left(\begin{array}{l}
z_{+}
\\
z_{-}
\end{array}\right)=z_{+}(0)\e^{xA}\left(\begin{array}{r}
1
\\
-1
\end{array}\right)+\int_{0}^{x}\e^{(x-\tau)A}\left(\begin{array}{r}
-1
\\
0
\end{array}\right)\cos(\lambda(\tau-\alpha))\,\ud \tau.
$$
Using the fact that $v(x)=z_{+}(x)+z_{-}(x)$, we follow
\begin{eqnarray}\label{pwkv132}
v(x)&=&2iz_{+}(0)\sin(\lambda x)-\int_{0}^{x}\e^{i\lambda(x-\tau)}\cos(\lambda(\tau-\alpha))\,\ud\tau\nonumber
\\
&=&2iz_{+}(0)\sin(\lambda x)-\frac{x}{2}\e^{i\lambda(x-\alpha)}-\frac{1}{2\lambda}\e^{i\lambda\alpha}\sin(\lambda x),
\end{eqnarray}
and in particular, we have
$$
v(\alpha^{-})=-\frac{\alpha}{2}+\sin(\lambda\alpha)\left[2iz_{+}(0)-\frac{1}{2\lambda}\e^{-i\lambda\alpha}\right].
$$
Furthermore, since $i\lambda u(\alpha^{-})=f(\alpha^{-})+v(\alpha^{-})=v(\alpha^{-})$ then we find that
\begin{equation}\label{pwkv138}
u(\alpha^{-})=\frac{i\alpha}{2\lambda}-\frac{i}{\lambda}\sin(\lambda\alpha)\left[2iz_{+}(0)-\frac{1}{2\lambda}\e^{-i\lambda\alpha}\right].
\end{equation}
Similarly, using $u'(x)=z_{+}(x)-z_{-}(x)$ then we get
\begin{equation*}
\begin{split}
u'(x)=2z_{+}(0)\cos(\lambda x)-\int_{0}^{x}\e^{i\lambda (x-\tau)}\cos(\lambda(\tau-\alpha))\,\ud\tau
\end{split}
\end{equation*}
and this leads to
\begin{equation}\label{pwkv139}
u'(\alpha^{-})=-\frac{\alpha}{2}+2\cos(\lambda\alpha)z_{+}(0)-\frac{\sin(\lambda\alpha)}{2\lambda}\e^{-i\lambda\alpha}.
\end{equation}

For $x\in(\alpha,\beta)$, we have
\begin{equation}\label{pwkv124}
\left\{\begin{array}{ll}
i\lambda u-v=0&\text{in }(\alpha,\beta)
\\
i\lambda v-(1+ia\lambda)u''=0&\text{in }(\alpha,\beta).
\end{array}\right.
\end{equation}
The solution of~\eqref{pwkv124} is
$$
u(x)=C_{1}\e^{\omega x}+C_{2}\e^{-\omega x}
$$
where $C_{1}$ and $C_{2}$ are constant and
$$
\omega=\omega_{n}=\frac{\lambda}{(1+a^{2}\lambda^{2})^{\frac{1}{4}}}\left(\cos\left(\frac{\theta}{2}\right)+i\sin\left(\frac{\theta}{2}\right)\right)
$$
with
$$\displaystyle\cos(\theta)=\frac{-1}{(1+a^{2}\lambda^{2})^{\frac{1}{2}}}\longrightarrow 0\,\text{  and  }\,\sin(\theta)=\frac{a\lambda}{(1+a^{2}\lambda^{2})^{\frac{1}{2}}}\longrightarrow 1\,\text{ as } n\longrightarrow+\infty.$$
By the continuity conditions at $x=\alpha$, i.e.
\begin{equation*}
\left\{\begin{array}{l}
\displaystyle u(\alpha^{+})=u(\alpha^{-})
\\
\displaystyle(1+ia\lambda)u'(\alpha^{+})=u'(\alpha^{-}),
\end{array}\right.
\end{equation*}
we find that
\begin{equation}\label{pwkv130}
C_{1}=\frac{1}{2}\e^{-\omega\alpha}\left(u(\alpha^{-})+\frac{u'(\alpha^{-})}{\omega(1+ia\lambda)}\right)\text{ and }C_{2}=\frac{1}{2}\e^{\omega\alpha}\left(u(\alpha^{-})-\frac{u'(\alpha^{-})}{\omega(1+ia\lambda)}\right).
\end{equation}

For $x\in(\beta,l)$, we have
\begin{equation}\label{pwkv125}
\left\{\begin{array}{ll}
i\lambda u-v=0&\text{in }(\beta,l)
\\
i\lambda v-u''=0&\text{in }(\beta,l).
\end{array}\right.
\end{equation}
The solution of~\eqref{pwkv125} is given by
$$
u(x)=C_{3}\e^{i\lambda(x-l)}+C_{4}\e^{-i\lambda(x-l)}.
$$
By continuity conditions at $x=\beta$, i.e.
\begin{equation*}
\left\{\begin{array}{l}
u(\beta^{-})=u(\beta^{+})
\\
(1+ia\lambda)u'(\beta^{-})=u'(\beta^{+}),
\end{array}\right.
\end{equation*}
we find that
\begin{equation}\label{pwkv128}
C_{3}=\frac{1}{2}\e^{-i\lambda(\beta-l)}\left[C_{1}\e^{\omega\beta}\left(1+\frac{\omega}{i\lambda}(1+ia\lambda)\right)+C_{2}\e^{-\omega\beta}\left(1-\frac{\omega}{i\lambda}(1+ia\lambda)\right)\right]
\end{equation}
and
\begin{equation}\label{pwkv129}
C_{4}=\frac{1}{2}\e^{i\lambda(\beta-l)}\left[C_{1}\e^{\omega\beta}\left(1-\frac{\omega}{i\lambda}(1+ia\lambda)\right)+C_{2}\e^{-\omega\beta}\left(1+\frac{\omega}{i\lambda}(1+ia\lambda)\right)\right].
\end{equation}

For $x\in(l,L)$, we have
\begin{equation}\label{pwkv126}
\left\{\begin{array}{ll}
i\lambda u-v=0&\text{in }(l,L)
\\
i\lambda v+u''''=0&\text{in }(l,L).
\end{array}\right.
\end{equation}
The solution of~\eqref{pwkv126} is given by
$$
u(x)=C_{5}\e^{\sqrt{\lambda}(x-l)}+C_{6}\e^{-\sqrt{\lambda}(x-l)}+C_{7}\e^{i\sqrt{\lambda}(x-l)}+C_{8}\e^{-i\sqrt{\lambda}(x-l)}
$$
where $C_{5},\,C_{6},\,C_{7}$ and $C_{8}$ by transmission and boundary conditions
\begin{equation*}
\left\{\begin{array}{l}
u(L)=u'(L)=u'(l^{+})=0
\\
 u(l^{-})=u(l^{+})
\\
u'(l^{-})=u'''(l^{+})
\end{array}\right.
\end{equation*}
satisfies
\begin{equation}\label{pwkv127}
\left\{\begin{array}{l}
C_{5}\e^{2n\pi}+C_{6}\e^{-2n\pi}+C_{7}+C_{8}=0
\\
C_{5}\e^{2n\pi}-C_{6}\e^{-2n\pi}+iC_{7}-iC_{8}=0
\\
C_{5}-C_{6}+iC_{7}-iC_{8}=0
\\
C_{5}+C_{6}+C_{7}+C_{8}=C_{3}+C_{4}
\\
C_{5}-C_{6}-iC_{7}+iC_{8}=\frac{i}{\sqrt{\lambda}}(C_{4}-C_{3}).
\end{array}\right.
\end{equation}
Consider only the four first equations of~\eqref{pwkv127} then we obtain
$$
C_{5}=\frac{(C_{3}+C_{4})(\e^{-2n\pi}-1)}{4(1-\cosh(2n\pi))},\qquad C_{6}=\frac{(C_{3}+C_{4})(\e^{2n\pi}-1)}{4(1-\cosh(2n\pi))},
$$
$$
C_{7}=\frac{(C_{3}+C_{4})(\cosh(2n\pi)-1-i\sinh(2n\pi))}{4(1-\cosh(2n\pi))},
$$
and
$$
C_{8}=\frac{(C_{3}+C_{4})(\cosh(2n\pi)-1+i\sinh(2n\pi))}{4(1-\cosh(2n\pi))}.
$$
We substitute all that terms into the fifth equation of~\eqref{pwkv127} and using the expressions of $C_{3}$ and $C_{4}$ in~\eqref{pwkv128} and~\eqref{pwkv129} we obtain
$$
\frac{\sinh(2n\pi)(C_{1}\e^{\omega\beta}+C_{2}\e^{-\omega\beta})}{2i(1-\cosh(2n\pi))}=\frac{\omega}{\lambda^{\frac{3}{2}}(1+ia\lambda)}(C_{1}\e^{\omega\beta}-C_{2}\e^{-\omega\beta}).
$$
Referring to the expression of $C_{1}$ and $C_{2}$ in~\eqref{pwkv130} we get
\begin{equation}\label{pwkv131}
\begin{split}
\left[\frac{i\omega\sinh(2n\pi)}{2(1-\cosh(2n\pi))}+\frac{\omega^{2}\tanh(\omega(\beta-\alpha))}{(1+ia\lambda)\lambda^{\frac{3}{2}}}\right]\left[\frac{i\alpha}{2}-i\sin(\lambda\alpha)\left(2iz_{+}(0)-\frac{1}{2\lambda}\e^{-i\lambda\alpha}\right)\right]=
\\
\left[\frac{\lambda\sinh(2n\pi)\tanh(\omega(\beta-\alpha))}{2i(1+ia\lambda)(1-\cosh(2n\pi))}+\frac{\omega^{2}}{\lambda^{\frac{1}{2}}(1+ia\lambda)^{2}}\right]\left[-\frac{\alpha}{2}+2\cos(\lambda\alpha)z_{+}(0)-\frac{\sin(\lambda\alpha)}{2\lambda}\e^{-i\lambda\alpha}\right],
\end{split}
\end{equation}
where here we have used the expressions of $u(\alpha^{-})$ and $u'(\alpha^{-})$ in~\eqref{pwkv138} and~\eqref{pwkv139} respectively. Since the left hand side of~\eqref{pwkv138} is equivalent to $\displaystyle i\omega\left(\frac{i\alpha}{2}+2z_{+}(0)\sin(\lambda\alpha)\right)$ and the right hand side is equivalent to $\displaystyle\frac{1}{2a}\left(2z_{+}(0)\cos(\lambda\alpha)-\frac{\alpha}{2}\right)$, using the fact that $\omega$ is equivalent to $\displaystyle i\sqrt{\frac{\lambda}{a}}$ we obtain
\begin{equation}\label{pwkv140}
\sqrt{\lambda}\left(\frac{i\alpha}{2}+2z_{+}(0)\sin(\lambda\alpha)\right)\sim\frac{1}{2\sqrt{a}}\left(\frac{\alpha}{2}-2z_{+}(0)\cos(\lambda\alpha)\right).
\end{equation}
In what's follow we want to prove that $|z_{+}(0)|\longrightarrow+\infty$. It's clear that if $\displaystyle\frac{i\alpha}{2}+2z_{+}(0)\sin(\lambda\alpha)$ don't tends to zero then $|z_{+}(0)|$ tends to $+\infty$. We suppose now that $\displaystyle\frac{i\alpha}{2}+2z_{+}(0)\sin(\lambda\alpha)$ converge to zero. If $\sin(\lambda\alpha)$ tends to zero then $|z_{+}(0)|$ tends to $+\infty$. If $\sin(\lambda\alpha)$ and $|z_{+}(0)|$ converge to non zero real numbers then $\re(z_{+}(0))$ should tends to zero and this is absurd according to the right hand side term of~\eqref{pwkv140} which make $|z_{+}(0)|\longrightarrow+\infty$ as $n\longrightarrow+\infty$.
\\
Referring to~\eqref{pwkv132} this leads to
\begin{equation}
\begin{split}
\int_{0}^{\alpha}|v(x)|^{2}\ud x&\geq\left(2|z_{+}(0)|^{2}-\frac{1}{4\lambda^{2}}\right)\int_{0}^{\alpha}\sin^{2}\left(\lambda x\right)\ud x-\frac{\alpha^{3}}{12}
\\
&\geq\left(2|z_{+}(0)|^{2}-\frac{1}{4\lambda^{2}}\right)\left(\frac{\alpha}{2}-\frac{\sin(2\lambda\alpha)}{4\lambda}\right)-\frac{\alpha^{3}}{12}\longrightarrow+\infty.
\end{split}
\end{equation}
Since $\displaystyle\|(f_{n},g_{n})\|_{\mathcal{H}_{2}}^{2}=\alpha+\frac{1}{\lambda}\sin(2\alpha\lambda)\longrightarrow\alpha$ as $n\longrightarrow+\infty$ and
$$
\|(i\lambda_{n}-\mathcal{A}_{2})^{-1}(f_{n},g_{n})\|_{\mathcal{H}_{2}}^{2}=\|(u_{n},v_{n})\|_{\mathcal{H}_{2}}^{2}\geq\int_{0}^{\alpha}|v_{1,n}(x)|^{2}\ud x\longrightarrow+\infty \text{ as } n\longrightarrow+\infty,
$$
we conclude that $\displaystyle\sup_{\lambda\in\R}\|(i\lambda-\mathcal{A}_{2})^{-1}\|=+\infty$, thus $e^{\mathcal{A}_{2}t}$ is not exponentially stable.
\subsection{Case when the damping arising from the beam equation}
Let $\displaystyle\lambda=\lambda_{n}=\frac{2n\pi}{(L-l)}$, $n=1,2,\dots,$. We take $b$ constant in $(\alpha,\beta)$ and
\begin{equation*}
f(x)=f_{n}(x)=\left\{\begin{array}{cl}
0&0<x<l
\\
\displaystyle\frac{1}{\lambda}\sin(\lambda(x-\alpha))&l<x<L
\end{array}\right.
\end{equation*}
and
\begin{equation*}
g(x)=g_{n}(x)=\left\{\begin{array}{cl}
0&0<x<l
\\
\displaystyle\cos(\lambda(x-\alpha))&l<x<L.
\end{array}\right.
\end{equation*}
We set $f_{1,n},\,f_{2,n}$ the restriction of $f_{n}$ over the intervals $(0,l)$ and $(l,L)$ respectively. Similarly $g_{1,n},\,g_{2,n}$ are the restriction of $g_{n}$ over the intervals $(0,l)$ and $(l,L)$ respectively. In the intervals $(0,\alpha),\,(\alpha,\beta),\,(\beta,l)$ and $(l,L)$ we solve the resolvent equation
$$
(i\lambda_{n}-\mathcal{A}_{1})\left(\begin{array}{l}
w_{1,n}
\\
w_{2,n}
\\
v_{1,n}
\\
v_{2,n}
\end{array}\right)=\left(\begin{array}{l}
f_{1,n}
\\
f_{2,n}
\\
g_{1,n}
\\
g_{2,n}
\end{array}\right)
$$
where $(w,v)=(w_{n},v_{n})=((w_{1,n},w_{2,n}),(v_{1,n},v_{2,n}))\in\mathcal{D}(\mathcal{A}_{1})$.

For $x\in(0,\alpha)$, we have
\begin{equation}\label{pwkv133}
\left\{\begin{array}{ll}
i\lambda w-v=0&\text{in }(0,\alpha)
\\
i\lambda v+w''''=0&\text{in }(0,\alpha)
\\
w(0)=w'(0)=0.
\end{array}\right.
\end{equation}
The solution of~\eqref{pwkv133} is
$$
w(x)=K_{1}\e^{\sqrt{\lambda}x}+K_{2}\e^{-\sqrt{\lambda}x}+K_{3}\e^{i\sqrt{\lambda}x}+K_{4}\e^{-i\sqrt{\lambda}x},
$$
where $K_{1},\,K_{2},\,K_{3}$ and $K_{4}$ verifies thanks to the boundary conditions
\begin{equation}\label{pwkv134}
K_{1}=\frac{1}{2}((1+i)K_{3}+(1-i)K_{4})\;\text{ and }\;K_{2}=\frac{1}{2}((1-i)K_{3}+(1+i)K_{4})
\end{equation}

For $x\in(\alpha,\beta)$, we have
\begin{equation}\label{pwkv135}
\left\{\begin{array}{ll}
i\lambda w-v=0&\text{in }(\alpha,\beta)
\\
i\lambda v+w''''+bv''''=0&\text{in }(\alpha,\beta).
\end{array}\right.
\end{equation}
The solution of~\eqref{pwkv134} is
$$
w(x)=K_{5}\e^{\omega(x-\alpha)}+K_{6}\e^{-\omega(x-\alpha)}+K_{7}\e^{i\omega(x-\alpha)}+K_{8}\e^{-i\omega(x-\alpha)},
$$
where 
$$\displaystyle\omega=\omega_{n}=\frac{\sqrt{\lambda}}{\sqrt{1+b^{2}\lambda^{2}}}\left(\cos\left(\frac{\theta}{4}\right)+i\sin\left(\frac{\theta}{4}\right)\right)\longrightarrow0\;\text{ as }\;n\longrightarrow+\infty$$ with
$$
\cos(\theta)=\frac{1}{\sqrt{1+b^{2}\lambda^{2}}}\;\;\text{ and }\;\;\sin(\theta)=\frac{-b\lambda}{\sqrt{1+b^{2}\lambda^{2}}}.
$$
By the continuity condition at $x=\alpha$ i.e.
\begin{equation*}
\left\{\begin{array}{l}
w(\alpha^{-})=w(\alpha^{+})
\\
w'(\alpha^{-})=w'(\alpha^{+})
\\
w''(\alpha^{-})=(1+ib\lambda)w''(\alpha^{+})
\\
w'''(\alpha^{-})=(1+ib\lambda)w'''(\alpha^{+}),
\end{array}\right.
\end{equation*}
we find
\begin{eqnarray*}
K_{5}&=&\frac{1}{4}\left(\frac{1+i}{2}u_{1}\e^{\sqrt{\lambda}\alpha}+\frac{1-i}{2}u_{2}\e^{-\sqrt{\lambda}\alpha}+u_{3}\e^{i\sqrt{\lambda}\alpha}\right)K_{3}
\\
&+&\frac{1}{4}\left(\frac{1-i}{2}u_{1}\e^{\sqrt{\lambda}\alpha}+\frac{1+i}{2}u_{2}\e^{-\sqrt{\lambda}\alpha}+u_{4}\e^{-i\sqrt{\lambda}\alpha}\right)K_{4}=v_{1}K_{3}+v_{2}K_{4},
\end{eqnarray*}
\begin{eqnarray*}
K_{6}&=&\frac{1}{4}\left(\frac{1+i}{2}u_{2}\e^{\sqrt{\lambda}\alpha}+\frac{1-i}{2}u_{1}\e^{-\sqrt{\lambda}\alpha}+u_{4}\e^{i\sqrt{\lambda}\alpha}\right)K_{3}
\\
&+&\frac{1}{4}\left(\frac{1-i}{2}u_{2}\e^{\sqrt{\lambda}\alpha}+\frac{1+i}{2}u_{1}\e^{-\sqrt{\lambda}\alpha}+u_{3}\e^{-i\sqrt{\lambda}\alpha}\right)K_{4}=v_{3}K_{3}+v_{4}K_{4},
\end{eqnarray*}
\begin{eqnarray*}
K_{7}&=&\frac{1}{4}\left(\frac{1+i}{2}u_{4}\e^{\sqrt{\lambda}\alpha}+\frac{1-i}{2}u_{3}\e^{-\sqrt{\lambda}\alpha}+u_{1}\e^{i\sqrt{\lambda}\alpha}\right)K_{3}
\\
&+&\frac{1}{4}\left(\frac{1-i}{2}u_{4}\e^{\sqrt{\lambda}\alpha}+\frac{1+i}{2}u_{3}\e^{-\sqrt{\lambda}\alpha}+u_{2}\e^{-i\sqrt{\lambda}\alpha}\right)K_{4}=v_{5}K_{3}+v_{6}K_{4},
\end{eqnarray*}
\begin{eqnarray*}
K_{8}&=&\frac{1}{4}\left(\frac{1+i}{2}u_{3}\e^{\sqrt{\lambda}\alpha}+\frac{1-i}{2}u_{4}\e^{-\sqrt{\lambda}\alpha}+u_{2}\e^{i\sqrt{\lambda}\alpha}\right)K_{3}
\\
&+&\frac{1}{4}\left(\frac{1-i}{2}u_{3}\e^{\sqrt{\lambda}\alpha}+\frac{1+i}{2}u_{4}\e^{-\sqrt{\lambda}\alpha}+u_{1}\e^{-i\sqrt{\lambda}\alpha}\right)K_{4}=v_{7}K_{3}+v_{8}K_{4},
\end{eqnarray*}
where
\begin{eqnarray*}
u_{1}&=&1+\frac{\sqrt{\lambda}}{\omega}+(1+ib\lambda)\left(\frac{\sqrt{\lambda}}{\omega}\right)^{2}+(1+ib\lambda)\left(\frac{\sqrt{\lambda}}{\omega}\right)^{3},
\\
u_{2}&=&1-\frac{\sqrt{\lambda}}{\omega}+(1+ib\lambda)\left(\frac{\sqrt{\lambda}}{\omega}\right)^{2}-(1+ib\lambda)\left(\frac{\sqrt{\lambda}}{\omega}\right)^{3},
\\
u_{3}&=&1+i\frac{\sqrt{\lambda}}{\omega}-(1+ib\lambda)\left(\frac{\sqrt{\lambda}}{\omega}\right)^{2}-i(1+ib\lambda)\left(\frac{\sqrt{\lambda}}{\omega}\right)^{3},
\\
u_{4}&=&1-i\frac{\sqrt{\lambda}}{\omega}-(1+ib\lambda)\left(\frac{\sqrt{\lambda}}{\omega}\right)^{2}+i(1+ib\lambda)\left(\frac{\sqrt{\lambda}}{\omega}\right)^{3}.
\end{eqnarray*}

For $x\in(\beta,l)$, we have
\begin{equation}\label{pwkv136}
\left\{\begin{array}{ll}
i\lambda w-v=0&\text{in }(\beta,l)
\\
i\lambda v+w''''=0&\text{in }(\beta,l)
\\
w'(l)=0.
\end{array}\right.
\end{equation}
The solution of~\eqref{pwkv136} is
$$
w(x)=K_{9}\e^{\sqrt{\lambda}(x-\beta)}+K_{10}\e^{-\sqrt{\lambda}(x-\beta)}+K_{11}\e^{i\sqrt{\lambda}(x-\beta)}+K_{12}\e^{-i\sqrt{\lambda}(x-\beta)},
$$
where 
$$\displaystyle\omega=\omega_{n}=\frac{\sqrt{\lambda}}{\sqrt{1+b^{2}\lambda^{2}}}\left(\cos\left(\frac{\theta}{4}\right)+i\sin\left(\frac{\theta}{4}\right)\right)\longrightarrow0\;\text{ as }\;n\longrightarrow+\infty$$ with
$$
\cos(\theta)=\frac{1}{\sqrt{1+b^{2}\lambda^{2}}}\;\;\text{ and }\;\;\sin(\theta)=\frac{-b\lambda}{\sqrt{1+b^{2}\lambda^{2}}}.
$$
By the continuity condition at $x=\beta$, i.e.
\begin{equation*}
\left\{\begin{array}{l}
w(\beta^{+})=w(\beta^{-})
\\
w'(\beta^{+})=w'(\beta^{-})
\\
w''(\beta^{+})=(1+ib\lambda)w''(\beta^{-})
\\
w'''(\beta^{+})=(1+ib\lambda)w'''(\beta^{-}),
\end{array}\right.
\end{equation*}
we find
\begin{eqnarray*}
K_{9}\!\!\!&=&\!\!\!\!\frac{1}{16}\left(v_{1}w_{1}\e^{\omega(\beta-\alpha)}+v_{3}w_{2}\e^{-\omega(\beta-\alpha)}+v_{5}w_{3}\e^{i\omega(\beta-\alpha)}+v_{7}w_{4}\e^{-i\omega(\beta-\alpha)}\right)K_{3}
\\
&+&\!\!\!\!\frac{1}{16}\left(v_{2}w_{1}\e^{\omega(\beta-\alpha)}+v_{4}w_{2}\e^{-\omega(\beta-\alpha)}+v_{6}w_{3}\e^{i\omega(\beta-\alpha)}+v_{8}w_{4}\e^{-i\omega(\beta-\alpha)}\right)K_{4}=z_{1}K_{3}+z_{2}K_{4},
\end{eqnarray*}
\begin{eqnarray*}
K_{10}\!\!\!&=&\!\!\!\!\frac{1}{16}\left(v_{1}w_{2}\e^{\omega(\beta-\alpha)}+v_{3}w_{1}\e^{-\omega(\beta-\alpha)}+v_{5}w_{4}\e^{i\omega(\beta-\alpha)}+v_{7}w_{3}\e^{-i\omega(\beta-\alpha)}\right)K_{3}
\\
&+&\!\!\!\!\frac{1}{16}\left(v_{2}w_{2}\e^{\omega(\beta-\alpha)}+v_{4}w_{1}\e^{-\omega(\beta-\alpha)}+v_{6}w_{4}\e^{i\omega(\beta-\alpha)}+v_{8}w_{3}\e^{-i\omega(\beta-\alpha)}\right)K_{4}=z_{3}K_{3}+z_{4}K_{4},
\end{eqnarray*}
\begin{eqnarray*}
K_{11}\!\!\!&=&\!\!\!\!\frac{1}{16}\left(v_{1}w_{4}\e^{\omega(\beta-\alpha)}+v_{3}w_{3}\e^{-\omega(\beta-\alpha)}+v_{5}w_{1}\e^{i\omega(\beta-\alpha)}+v_{7}w_{2}\e^{-i\omega(\beta-\alpha)}\right)K_{3}
\\
&+&\!\!\!\!\frac{1}{16}\left(v_{2}w_{4}\e^{\omega(\beta-\alpha)}+v_{4}w_{3}\e^{-\omega(\beta-\alpha)}+v_{6}w_{1}\e^{i\omega(\beta-\alpha)}+v_{8}w_{2}\e^{-i\omega(\beta-\alpha)}\right)K_{4}=z_{5}K_{3}+z_{6}K_{4},
\end{eqnarray*}
\begin{eqnarray*}
K_{12}\!\!\!&=&\!\!\!\!\frac{1}{16}\left(v_{1}w_{3}\e^{\omega(\beta-\alpha)}+v_{3}w_{4}\e^{-\omega(\beta-\alpha)}+v_{5}w_{2}\e^{i\omega(\beta-\alpha)}+v_{7}w_{1}\e^{-i\omega(\beta-\alpha)}\right)K_{3}
\\
&+&\!\!\!\!\frac{1}{16}\left(v_{2}w_{3}\e^{\omega(\beta-\alpha)}+v_{4}w_{4}\e^{-\omega(\beta-\alpha)}+v_{6}w_{2}\e^{i\omega(\beta-\alpha)}+v_{8}w_{1}\e^{-i\omega(\beta-\alpha)}\right)K_{4}=z_{7}K_{3}+z_{8}K_{4},
\end{eqnarray*}
where
\begin{eqnarray*}
w_{1}&=&1+\frac{\omega}{\sqrt{\lambda}}+(1+ib\lambda)\left(\frac{\omega}{\sqrt{\lambda}}\right)^{2}+(1+ib\lambda)\left(\frac{\omega}{\sqrt{\lambda}}\right)^{3}\longrightarrow 1,
\\
w_{2}&=&1-\frac{\omega}{\sqrt{\lambda}}+(1+ib\lambda)\left(\frac{\omega}{\sqrt{\lambda}}\right)^{2}-(1+ib\lambda)\left(\frac{\omega}{\sqrt{\lambda}}\right)^{3}\longrightarrow 1,
\\
w_{3}&=&1+i\frac{\omega}{\sqrt{\lambda}}-(1+ib\lambda)\left(\frac{\omega}{\sqrt{\lambda}}\right)^{2}-i(1+ib\lambda)\left(\frac{\omega}{\sqrt{\lambda}}\right)^{3}\longrightarrow 1,
\\
w_{4}&=&1-i\frac{\omega}{\sqrt{\lambda}}-(1+ib\lambda)\left(\frac{\omega}{\sqrt{\lambda}}\right)^{2}+i(1+ib\lambda)\left(\frac{\omega}{\sqrt{\lambda}}\right)^{3}\longrightarrow 1.
\end{eqnarray*}
By taking account the condition $w'(l^{-})=0$ we obtain
\begin{equation*}
\begin{split}
w(x)=\frac{K_{4}}{16}\bigg[Z\times\left(z_{1}\e^{\sqrt{\lambda}(x-\beta)}+z_{3}\e^{-\sqrt{\lambda}(x-\beta)}+z_{5}\e^{i\sqrt{\lambda}(x-\beta)}+z_{7}\e^{-i\sqrt{\lambda}(x-\beta)}\right)
\\
+z_{2}\e^{\sqrt{\lambda}(x-\beta)}+z_{4}\e^{-\sqrt{\lambda}(x-\beta)}+z_{6}\e^{i\sqrt{\lambda}(x-\beta)}+z_{8}\e^{-i\sqrt{\lambda}(x-\beta)}\bigg],
\end{split}
\end{equation*}
where
$$
Z=\frac{-z_{2}\e^{\sqrt{\lambda}(l-\beta)}+z_{4}\e^{-\sqrt{\lambda}(l-\beta)}-iz_{6}\e^{i\sqrt{\lambda}(l-\beta)}+iz_{8}\e^{-i\sqrt{\lambda}(l-\beta)}}{z_{1}\e^{\sqrt{\lambda}(l-\beta)}-z_{3}\e^{-\sqrt{\lambda}(l-\beta)}+iz_{5}\e^{i\sqrt{\lambda}(l-\beta)}-iz_{7}\e^{-i\sqrt{\lambda}(l-\beta)}}.
$$

For $x\in(l,L)$, we have
\begin{equation}\label{pwkv137}
\left\{\begin{array}{ll}
i\lambda w-v=f&\text{in }(l,L)
\\
i\lambda v-w''=g&\text{in }(l,L)
\\
w(L)=0.
\end{array}\right.
\end{equation}
By proceeding as in the previous subsection (case when $x\in(0,\alpha)$) then we find
$$
w(l^{+})=\frac{L-l}{2i\lambda}\;\text{ and }\;w'(l^{+})=2z_{+}(L)+\frac{L-l}{2},
$$
where
$$
z_{+}(x)=\frac{1}{2}\left(v(x)+w'(x)\right).
$$
The transmission conditions
\begin{equation*}
\left\{\begin{array}{l}
w(l^{+})=w(l^{-})
\\
w'(l^{+})=-w'''(l^{-})
\end{array}\right.
\end{equation*}
allows to solve $K_{4}$ and particularly to obtain
\begin{equation*}
\begin{split}
2z_{+}(L)+\frac{L-l}{2}=-\frac{\sqrt{\lambda}(L-l)}{2i}\Big[Z\left(z_{1}\e^{\sqrt{\lambda}(l-\beta)}+z_{3}\e^{-\sqrt{\lambda}(l-\beta)}+z_{5}\e^{i\sqrt{\lambda}(l-\beta)}+z_{7}\e^{-i\sqrt{\lambda}(l-\beta)}\right)
\\
+z_{2}\e^{\sqrt{\lambda}(l-\beta)}+z_{4}\e^{-\sqrt{\lambda}(l-\beta)}+z_{6}\e^{i\sqrt{\lambda}(l-\beta)}+z_{8}\e^{-i\sqrt{\lambda}(l-\beta)}\Big]^{-1}\times\Big[Z\Big(z_{1}\e^{\sqrt{\lambda}(l-\beta)}-z_{3}\e^{-\sqrt{\lambda}(l-\beta)}
\\
-iz_{5}\e^{i\sqrt{\lambda}(l-\beta)}+iz_{7}\e^{-i\sqrt{\lambda}(l-\beta)}\Big)+z_{2}\e^{\sqrt{\lambda}(l-\beta)}-z_{4}\e^{-\sqrt{\lambda}(l-\beta)}-iz_{6}\e^{i\sqrt{\lambda}(l-\beta)}+iz_{8}\e^{-i\sqrt{\lambda}(l-\beta)}\Big].
\end{split}
\end{equation*}
This further leads since $\lambda\longrightarrow+\infty$ to
$$
|z_{+}(L)|\longrightarrow+\infty\text{ as }n\longrightarrow+\infty.
$$
Using now the same arguments as in the end of the previous subsection to obtain 
$$\displaystyle\sup_{\lambda\in\R}\|(i\lambda-\mathcal{A}_{1})^{-1}\|=+\infty$$
and hence $e^{\mathcal{A}_{1}t}$ is not exponentially stable.
\subsubsection*{Acknowledgments}
The author thanks the referees for many valuable remarks which helped us to improve the paper significantly.
\nocite{*}
\bibliographystyle{alpha}
\bibliography{BibPWKV1D}
\addcontentsline{toc}{section}{References}
\end{document}